\documentclass[a4paper,10pt,reqno]{amsart}
\usepackage{amsmath,amsfonts,amscd,amssymb,graphicx,mathrsfs,eufrak,dsfont}
\usepackage[dvips,all,arc,curve,color,frame]{xy}
\usepackage[usenames]{color}

\usepackage[colorlinks]{hyperref}
\usepackage{tikz,mathrsfs}
\usetikzlibrary{arrows,decorations.pathmorphing,decorations.pathreplacing,positioning,shapes.geometric,shapes.misc,decorations.markings,decorations.fractals,calc,patterns}

\tikzset{>=stealth',
       cvertex/.style={circle,draw=black,inner sep=1pt,outer sep=3pt},
       vertex/.style={circle,fill=black,inner sep=1pt,outer sep=3pt},
       star/.style={circle,fill=yellow,inner sep=0.75pt,outer sep=0.75pt},
       tvertex/.style={inner sep=1pt,font=\scriptsize},
       gap/.style={inner sep=0.5pt,fill=white}}

\newcommand{\arrowrl}[3][20]
{
\hspace{-5pt}
\begin{tikzpicture}
 \node (A) at (0,0) {};
 \node (B) at (1,0) {};
 \draw[->] ($(A)+(0,0.2)$) -- node [above] {$\scriptstyle f^*$} ($(B)+(0,0.2)$);
\draw [->] ($(B)+(0,0.2)$) -- node [below] {$\scriptstyle f_*$} ($(A)+(0,0.2)$);
\end{tikzpicture}
\hspace{-5pt}
}
\newcommand{\adj}[2][20]{\arrowrl}

\newtheorem{thm}{Theorem}[section]
\newtheorem{prop}[thm]{Proposition}
\newtheorem{lemma}[thm]{Lemma}
\newtheorem{defin}[thm]{Definition}
\newtheorem{cor}[thm]{Corollary}

\theoremstyle{definition} 
 
\newtheorem{example}[thm]{Example}

\newtheorem{remark}[thm]{Remark}
\newtheorem{conj}[thm]{Conjecture}

\newcommand{\F}{\mathcal{F}}
\newcommand{\G}{\mathcal{G}}
\renewcommand{\O}{\mathcal{O}}

\newcommand{\C}[1]{\mathbb{C}^{#1}}
\newcommand{\m}{\mathfrak{m}}
\newcommand{\n}{\mathfrak{n}}
\newcommand{\p}{\mathfrak{p}}

\renewcommand{\c}[1]{\mathcal{#1}}

\newcommand{\A}{\mathscr{A}}

\renewcommand{\t}[1]{\textnormal{#1}}

\newcommand{\looptop}[2]{\xy \SelectTips{cm}{10}
\POS(0,0) \endxy}

\newcommand{\thick}{\mathop{\rm thick}\nolimits}

\def\RHom{\mathop{\rm {\bf R}Hom}\nolimits}

\def\op{\mathop{\rm op}\nolimits}
\def\GL{\mathop{\rm GL}\nolimits}
\def\CM{\mathop{\rm CM}\nolimits}
\def\uCM{\mathop{\underline{\rm CM}}\nolimits}

\def\depth{\mathop{\rm depth}\nolimits}
\def\fl{\mathop{\sf fl}\nolimits}
\def\hgt{\mathop{\rm ht}\nolimits}
\def\mod{\mathop{\rm mod}\nolimits}

\def\coh{\mathop{\rm coh}\nolimits}
\def\Mod{\mathop{\rm Mod}\nolimits}
\def\Qcoh{\mathop{\rm Qcoh}\nolimits}
\def\refl{\mathop{\rm ref}\nolimits}
\def\proj{\mathop{\rm proj}\nolimits}
\def\Proj{\mathop{\rm Proj}\nolimits}
\def\FPD{\mathop{\rm FPD}\nolimits}
\def\Free{\mathop{\rm Free}\nolimits}

\def\pd{\mathop{\rm proj.dim}\nolimits}
\def\fd{\mathop{\rm flat.dim}\nolimits}
\def\id{\mathop{\rm inj.dim}\nolimits}
\def\Hom{\mathop{\rm Hom}\nolimits}
\def\End{\mathop{\rm End}\nolimits}
\def\Ext{\mathop{\rm Ext}\nolimits}

\def\add{\mathop{\rm add}\nolimits}

\def\Ker{\mathop{\rm Ker}\nolimits}
\def\Im{\mathop{\rm Im}\nolimits}

\def\Supp{\mathop{\rm Supp}\nolimits}

\def\Ann{\mathop{\rm Ann}\nolimits}
\def\Spec{\mathop{\rm Spec}\nolimits}
\def\Max{\mathop{\rm Max}\nolimits}
\def\gl{\mathop{\rm gl.dim}\nolimits}

\def\D{\mathop{\rm{D}^{}}\nolimits}
\def\Db{\mathop{\rm{D}^b}\nolimits}
\def\Dbf{\mathop{\rm{D}^b_{\fl}}\nolimits}
\def\Dsg{\mathop{\rm{D}_{\mathrm{sg}}}\nolimits}
\def\DsgQ{\mathop{\rm{D}_{\mathsf{SG}}}\nolimits}
\def\Kb{\mathop{\rm{K}^b}\nolimits}
\def\Perf{\mathop{\rm{per}}\nolimits}
\def\Lfr{\mathop{\mathrm{Lfr}}\nolimits}

\def\Rf{\mathop{\mathbf{R}f_*}\nolimits}
\def\Rg{\mathop{\mathbf{R}g_*}\nolimits}

\def\Lf{\mathop{{\bf L}f^*}\nolimits}

\def\fsh{\mathop{f^!}\nolimits}
\def\gsh{\mathop{g^!}\nolimits}
\def\RHom{\mathop{\rm {\bf R}Hom}\nolimits}

\def\RsHom{\mathop{{\bf R}\mathcal{H}om}\nolimits}
\def\RG{\mathop{{\bf R}\Gamma}\nolimits}

\newcommand{\cd}{{\mathcal D}}
\newcommand{\B}{\mathscr{B}}
\def\pretr{\mathop{\rm pretr}\nolimits}
\def\dg{\mathop{\rm dg}\nolimits}
\def\rC{\mathop{{}_{}\rm C}\nolimits}

\renewcommand{\k}{k}

\newcommand{\KK}{\mathbb{K}}
\newcommand{\kk}{\mathds{k}}

\edef\marginnotetextwidth{\the\textwidth}

\begin{document}
\title[\textsc{$\mathds{Q}$-factorial terminalizations and MMAs}]{\textsc{Singular Derived Categories of $\mathds{Q}$-factorial terminalizations and Maximal Modification Algebras}}
\author{Osamu Iyama}
\address{Osamu Iyama\\ Graduate School of Mathematics\\ Nagoya University\\ Chikusa-ku, Nagoya, 464-8602, Japan}
\email{iyama@math.nagoya-u.ac.jp}
\author{Michael Wemyss}
\address{Michael Wemyss, School of Mathematics, James Clerk Maxwell Building, The King's Buildings, Mayfield Road, Edinburgh, EH9 3JZ, UK.}
\email{wemyss.m@googlemail.com}
\begin{abstract}
Let $X$ be a Gorenstein normal 3-fold satisfying (ELF) with local rings which are at worst isolated hypersurface (e.g.\ terminal) singularities. By using the singular derived category $\Dsg(X)$ and its idempotent completion $\overline{\Dsg(X)}$, we give necessary and sufficient categorical conditions for $X$ to be $\mathds{Q}$-factorial and complete locally $\mathds{Q}$-factorial respectively.  We then relate this information to maximal modification algebras(=MMAs), introduced in \cite{IW}, by showing that if an algebra $\Lambda$ is derived equivalent to $X$ as above, then $X$ is $\mathds{Q}$-factorial if and only if $\Lambda$ is an MMA. Thus all rings derived equivalent to $\mathds{Q}$-factorial terminalizations in dimension three are MMAs.   As an application, we extend some of the algebraic results in \cite{BIKR} and \cite{DH} using geometric arguments.
\end{abstract}
\maketitle
\parindent 20pt
\parskip 0pt

\tableofcontents

\section{Introduction}

\subsection{Overview}\label{philosophy} 

The broader purpose of this work, in part a continuation of \cite{IW}, is to understand (and run) certain aspects of the minimal model program (=MMP) in dimension three using categorical and noncommutative techniques.  As part of the conjectural underlying picture, it is believed in dimension three (see \ref{3foldconj} below) that the theory of noncommutative minimal models (=MMAs below) should `control' the geometry of commutative minimal models (=$\mathds{Q}$-factorial terminalizations) in the same way noncommutative crepant resolutions (=NCCRs) control the geometry of crepant resolutions \cite{VdBNCCR}.  If this is to have any chance of being true, we must first  be able to extend the well-known characterization of smoothness by the singular derived category to be zero to a characterization of $\mathds{Q}$-factoriality.  This involves moving from smooth schemes to singular ones, and it is unclear geometrically what form such a derived category characterization should take.

However, from the noncommutative side, in the study of MMAs one homological condition, which we call `rigid-freeness', completely characterizes when `maximality' has been reached.  The purpose of this paper is to show that this very same condition gives a necessary condition for a scheme with isolated Gorenstein singularities to be $\mathds{Q}$-factorial. If furthermore the singularities are hypersurfaces (e.g.\ Gorenstein terminal singularities) this condition is both necessary and sufficient.  The ability of the noncommutative side to give new results purely in the setting of algebraic geometry adds weight to the conjectures outlined in \S\ref{conjectures} below.

We now explain our results in more detail.

\subsection{Results on $\mathds{Q}$-factoriality}\label{Qfactsection_intro}
We denote by $\CM R$ the category of CM (=maximal Cohen-Macaulay) $R$-modules, and by $\uCM R$ its stable category
(see \S\ref{prelim}). As in \cite{Orlov, OrlovCompletion}, we denote $\Dsg(X):=\Db(\coh X)/\Perf(X)$ (see \S\ref{Dsgprelim} for more details). We call a functor $F\colon \c{T}\to\c{T}'$ an \emph{equivalence up to direct summands} if it is fully faithful and any object $s\in\c{T}'$ is isomorphic to a direct summand of $F(t)$ for some $t\in\c{T}$. The first key technical result is the following.  Note that the latter part on the idempotent completion also follows from work of Orlov \cite{OrlovCompletion}, and is well-known to experts \cite{BK}.

\begin{thm}\label{Dsgsplits_intro}(=\ref{Dsgsplits})
Suppose that $X$ is a Gorenstein scheme of dimension $d$ satisfying (ELF) (see \ref{Dsgsplits} for full details and explanation), with isolated singularities $\{x_1,\hdots,x_n\}$.   Then there is a triangle equivalence
\[
\Dsg(X)\hookrightarrow\bigoplus_{i=1}^{n}\uCM\c{O}_{X,x_i}
\]
up to direct summands.  Thus taking the idempotent completion gives a triangle equivalence
\[
\overline{\Dsg(X)}\simeq\bigoplus_{i=1}^{n}\uCM\widehat{\c{O}}_{X,x_i} .
\]
\end{thm}
Recall that a normal scheme $X$ is defined to be {\em $\mathds{Q}$-factorial} if for every Weil divisor $D$, there exists $n\in\mathbb{N}$ for which $nD$ is Cartier.  If we can always take $n=1$, we say that $X$ is {\em locally factorial}.  Note that these conditions are Zariski local,  i.e.\ $X$ is $\mathds{Q}$-factorial (respectively, locally factorial) if and only if
$\c{O}_{X,x}$ is $\mathds{Q}$-factorial (respectively, factorial)
for all closed points $x\in X$. It is well--known (and problematic) that, unlike smoothness, $\mathds{Q}$-factoriality cannot in general be checked complete locally.  Thus we define $X$ to be {\em complete locally $\mathds{Q}$-factorial} (respectively, {\em complete locally factorial}) if the completion $\widehat{\O}_{X,x}$ is $\mathds{Q}$-factorial (respectively, factorial) for all closed points $x\in X$.

Now we recall a key notion from representation theory.  Let $\c{T}$ be a triangulated category with a suspension functor $[1]$.  An object $a\in \c{T}$ is called \emph{rigid} if $\Hom_\c{T}(a,a[1])=0$.  We say that $\c{T}$ is \emph{rigid-free} if every rigid object in $\c{T}$ is isomorphic to the zero object. Applying \ref{Dsgsplits_intro}, we have the following result.

\begin{thm}\label{main_intro}(=\ref{locallyfactorial}, \ref{Section3_main}) Suppose that $X$ is a normal 3-dimensional Gorenstein scheme over a field $k$, satisfying (ELF), with isolated singularities $\{x_1,\hdots,x_n\}$. \\
\t{(1)} If $\Dsg(X)$ is rigid-free, then $X$ is locally factorial.\\
\t{(2)} If $\overline{\Dsg(X)}$ is rigid-free, then $X$ is complete
locally factorial.\\
\t{(3)} If $\c{O}_{X,x_i}$ are hypersurfaces for all $1\leq i\leq n$,  then the following are equivalent:
\begin{enumerate}
\item[(a)] $X$ is locally factorial.
\item[(b)] $X$ is $\mathds{Q}$-factorial.
\item[(c)] $\Dsg(X)$ is rigid-free.
\item[(d)] $\uCM\O_{X,x}$ is rigid-free for all closed points $x\in X$.
\end{enumerate}
\t{(4)} If $\widehat{\c{O}}_{X,x_i}$ are hypersurfaces for all $1\leq i\leq n$,  then the following are equivalent: 
\begin{enumerate}
\item[(a)] $X$ is complete locally factorial.
\item[(b)] $X$ is complete locally $\mathds{Q}$-factorial.
\item[(c)] $\overline{\Dsg(X)}$ is rigid-free.
\item[(d)] $\uCM\widehat{\O}_{X,x}$ is rigid-free for all closed points $x\in X$.
\end{enumerate}
\end{thm}

Thus the technical geometric distinction between a $\mathds{Q}$-factorial and complete locally $\mathds{Q}$-factorial variety is explained categorically by the existence of a triangulated category in which there are no rigid objects, but whose idempotent completion has many.  Note also that \ref{main_intro} demonstrates that by passing to the idempotent completion of $\Dsg(X)$ we lose information on the global geometry (see e.g.\ \ref{ex:MMnotcompletelocal}), and so thus should be avoided.

\subsection{Results on MMAs}
The homological characterization of rigid-freeness in \ref{main_intro} originated in the study of modifying and maximal modifying modules (see \ref{MMdef} for definitions).   Recall that if $M$ is a maximal modifying $R$-module, we say that $\End_R(M)$ is a \emph{maximal modification algebra (=MMA)}.
\begin{prop}\label{MMtestforuCM_intro}(=\ref{MMtestforuCM}, \ref{Buchweitz})
Suppose that $R$ is normal, Gorenstein, equi-codimensional, 3-dimensional ring, and let $M$ be a modifying $R$-module. Assume that $\End_R(M)$ has only isolated singularities.  Then $\End_R(M)$ is an MMA if and only if the category $\Dsg(\End_R(M))$ is rigid-free.
\end{prop}

Throughout this paper we will use the word \emph{variety} to mean
\begin{itemize}
\item a normal integral scheme, of finite type over a field $k$, satisfying (ELF) (see \S\ref{Dsgprelim}).
\end{itemize}
Recall that if $X$ and $Y$ are varieties, then a projective birational morphism $f\colon Y\to X$ is called {\em crepant} if $f^*\omega_X=\omega_Y$ (see \S\ref{sect4} for more details). A {\em $\mathds{Q}$-factorial terminalization} of $X$ is a crepant projective birational morphism $f\colon Y\to X$ such that $Y$ has only $\mathds{Q}$-factorial terminal singularities. When $Y$ is furthermore smooth, we call $f$ a {\em crepant resolution}.  

Maximal modification algebras were introduced in \cite{IW} with the aim of generalizing NCCRs to cover the more general situation when crepant resolutions do not exist.  Throughout this paper, we say that $Y$ is {\em derived equivalent} to $\Lambda$ if $\Db(\coh Y)$ is equivalent to $\Db(\mod\Lambda)$.  The following result ensures that all algebras derived equivalent to a variety $Y$ above are of the form $\End_R(M)$, which relate the geometry to NCCRs and MMAs.

\begin{thm}\label{correct_form}(=\ref{Endincorrectform})
Let $Y\to\Spec R$ be a projective birational morphism between $d$-dimensional varieties.  Suppose that $Y$ is derived equivalent to some ring $\Lambda$.  If $\Lambda\in\refl R$, then $\Lambda\cong\End_R(M)$ for some $M\in\refl R$.
\end{thm}

Using this, we have the following result, which explains the geometric origin of the definition of modifying modules:

\begin{thm}\label{CMiif_intro}
Let $f\colon Y\to\Spec R$ be a projective birational morphism between $d$-dimensional Gorenstein varieties.  Suppose that $Y$ is derived equivalent to some ring $\Lambda$, then\\
\t{(1)} (=\ref{CMiffcrepant}) $f$ is crepant $\iff$ $\Lambda\in\CM R$.\\
\t{(2)} (=\ref{NCCRdimd}) $f$ is a crepant resolution $\iff$ $\Lambda$ is an NCCR of $R$.\\
In either case, $\Lambda\cong\End_R(M)$ for some $M\in\refl R$.
\end{thm}

We stress that \ref{CMiif_intro}(1) holds even if $Y$ is singular.   The following consequence of \ref{main_intro}, \ref{MMtestforuCM_intro} and \ref{CMiif_intro} is one of our main results.  

\begin{thm}\label{Qfact_intro}
(=\ref{mainQfact}) 
Let $f\colon Y\to\Spec R$ be a projective birational morphism, where $Y$ and $R$ are both Gorenstein varieties of dimension three. Assume that $Y$ has (at worst) isolated singularities $\{ x_1,\hdots,x_n\}$ where each $\c{O}_{Y,x_i}$ is a hypersurface.  If $Y$ is derived equivalent to some ring $\Lambda$, then the following are equivalent\\
\t{(1)} $f$ is crepant and $Y$ is $\mathds{Q}$-factorial.\\
\t{(2)} $\Lambda$ is an MMA of $R$.\\
In this situation, all MMAs of $R$ have isolated singularities, and are all derived equivalent.
\end{thm}

As an important consequence, when $\k=\mathbb{C}$ and $Y$ has only terminal singularities and derived equivalent to some ring $\Lambda$, then $Y$ is a $\mathds{Q}$-factorial terminalization of $\Spec R$ if and only if $\Lambda$ is an MMA.  See \ref{fieldC} for more details.

\subsection{Application to $cA_n$ singularities}
As an application of the above, we can extend some results of \cite{BIKR} and \cite{DH}.  We let $\KK$ denote an algebraically closed field of characteristic zero, and suppose that $f_1,\hdots,f_n\in\m:=(x,y)\subseteq \KK[[x,y]]$ are irreducible polynomials, and let
\[
R:=\KK[[u,v,x,y]]/(uv-f_1\hdots f_n).
\]
When $\KK=\mathbb{C}$ then $R$ is a $cA_m$ singularity for $m:={\rm ord}(f_1\hdots f_n)-1$ \cite[6.1(e)]{BIKR}. Note that $R$ is not assumed to be an isolated singularity, so $f_1,\hdots,f_n$ are not necessarily pairwise distinct. For each element $\omega$ in the symmetric group $\mathfrak{S}_n$, define
\[
T^\omega:=R\oplus (u,f_{\omega(1)})\oplus(u,f_{\omega(1)}f_{\omega(2)})\oplus\hdots\oplus (u,f_{\omega(1)}\hdots f_{\omega(n-1)}).
\]
Our next result generalizes \cite[1.5]{BIKR} and \cite[4.2]{DH}, and follows more or less immediately from \ref{Qfact_intro}.  This is stronger and more general since the results in \cite{BIKR, DH} assume that $R$ is an isolated singularity, and only study cluster tilting objects.

\begin{thm}\label{derivedMMCTforcAn_intro}(=\ref{derivedMMCTforcAn})
Let $f_1,\hdots,f_n\in\m:=(x,y)\subseteq \KK[[x,y]]$ be irreducible polynomials and $R=\KK[[u,v,x,y]]/(uv-f_1\hdots f_n)$. Then\\
\t{(1)}  Each $T^\omega$ ($\omega\in\mathfrak{S}_n$) is an MM $R$-module which is a generator.  The endomorphism rings $\End_R(T^\omega)$ have isolated singularities.\\
\t{(2)}  $T^\omega$ is a CT $R$-module for some $\omega\in\mathfrak{S}_n$ $\iff$ $T^\omega$ is a CT $R$-module for all $\omega\in\mathfrak{S}_n$ $\iff$ $f_i\notin\m^2$ for all $1\leq i\leq n$.
\end{thm}

We remark that the proof of \ref{derivedMMCTforcAn_intro} also provides a conceptual geometric reason for the condition $f_i\notin\m^2$.  We refer the reader to \S\ref{Dbsection} for details.

\subsection{Conjectures}\label{conjectures}

One of the motivations for the introduction of MMAs is the following conjecture, 
which naturally generalizes conjectures of Bondal--Orlov and Van den Bergh:
\begin{conj}\label{3foldconj}
Suppose that $R$ is a normal Gorenstein 3-fold over $\C{}$, with only canonical (equivalently rational, since $R$ is Gorenstein \cite[(3.8)]{YPG}) singularities.  Then\\
(1) $R$ admits an MMA.\\
(2) All $\mathds{Q}$-factorial terminalizations of $R$ and all MMAs of $R$ are derived equivalent.
\end{conj}
A special case of \ref{3foldconj} is the long--standing conjecture of Van den Bergh \cite[4.6]{VdBNCCR}, namely all crepant resolutions of $\Spec R$ (both commutative and non-commutative) are derived equivalent. The results in this paper (specifically \ref{Qfact_intro}) add some weight to \ref{3foldconj}.

We note that \ref{3foldconj} is true in some situations: 
\begin{thm}\label{evidence}
Suppose that $R$ is a normal Gorenstein 3-fold over $\mathbb{C}$ whose $\mathds{Q}$-factorial terminalizations $Y\to\Spec R$ have only one dimensional fibres (e.g.\ $R$ is a terminal Gorenstein singularity, or a $cA_{n}$ singularity in \S\ref{Dbsection}).  Then Conjecture \ref{3foldconj} is true.  
\end{thm}
\begin{proof}
The $\mathds{Q}$-factorial terminalizations are connected by a finite sequence of flops \cite{Kollar}, so they are derived equivalent by \cite{Chen}.  By \cite{VdB1d} there is a derived equivalence between $Y$ and some algebra $\End_R(M)$ with $\End_R(M)\in\CM R$.  Thus by \ref{Qfact_intro} $\End_R(M)$ is an MMA. On the other hand, we already know that MMAs are connected by tilting modules \cite[4.14]{IW}, so they are all derived equivalent. 
\end{proof}
As it stands, the proof of \ref{evidence} heavily uses the MMP.  We believe that it should be possible to prove \ref{evidence} (for more general situations) without using the MMP, instead building on the categorical techniques developed in this paper, and generalizing Van den Bergh's \cite{VdB1d} interpretation of Bridgeland--King--Reid \cite{BKR}.  Indeed, it is our belief and long-term goal to show that many of the results of the MMP come out of our categorical picture (along the lines of \cite{Bridgeland}, \cite{BKR}), allowing us to both bypass some of the classifications used in the MMP, and also run some aspects of the MMP in a more efficient manner.

\medskip
\noindent
{\bf Conventions.}  Throughout commutative rings are always assumed to be noetherian, and $R$ will always denote a commutative noetherian ring.  All modules will be left modules, so for a ring $A$  we denote $\mod A$ to be the category of finitely generated left $A$-modules, and $\Mod A$ will denote the category of all left $A$-modules.  Throughout when composing maps $fg$ will mean $f$ then $g$. Note that with these conventions $\Hom_R(M,X)$ is a $\End_R(M)$-module and $\Hom_R(X,M)$ is a $\End_R(M)^{\rm op}$-module.  For $M\in\mod A$ we denote $\add M$ to be the full subcategory consisting of direct summands of finite direct sums of copies of $M$, and we denote $\proj A:=\add A$ to be the category of finitely generated projective $A$-modules.  We say that $M\in\mod R$ is a {\em generator} if $R\in \add M$.  For an abelian category $\c{A}$, we denote by $\Db(\c{A})$ the bounded derived category of $\c{A}$.

If $X$ is a scheme, $\c{O}_{X,x}$ will denote the localization of the structure sheaf at the closed point $x\in X$.  We will denote by $\widehat{\c{O}}_{X,x}$ the completion of $\c{O}_{X,x}$ at the unique maximal ideal. For us, {\em locally} will always mean Zariski locally, that is if we say that $R$ has locally only isolated hypersurface singularities, we mean that each $R_\m$ is an isolated hypersurface singularity (or is smooth).  When we want to discuss the completion, we will always refer to this as {\em complete locally}.

Throughout, $\k$ will denote an arbitrary field.  Rings and schemes will {\em not} be assumed to be finite type over $\k$, unless specified.  When we say `over $\k$' we mean `of finite type over $\k$'.  

\medskip
\noindent
{\bf Acknowledgements.}   Thanks are due to Vanya Cheltsov, Anne-Sophie Kaloghiros, Raphael Rouquier and  Ed Segal for many invaluable suggestions and discussions.

\section{Preliminaries}\label{prelim}

For a commutative noetherian local ring $(R,\m)$ and $M\in\mod R$, recall that the \emph{depth} of $M$ is defined to be $\depth_R M:=\inf \{ i\geq 0: \Ext^i_R(R/\m,M)\neq 0 \}$.  We say that $M\in\mod R$ is \emph{maximal Cohen-Macaulay} (or simply, \emph{CM}) if $\depth_R M=\dim R$.  Now let $R$ be a (not necessarily local) commutative noetherian ring.  We say that $M\in\mod R$ is \emph{CM} if $M_\p$ is CM for all prime ideals $\p$ in $R$, and we say that $R$ is a \emph{CM ring} if $R$ is a CM $R$-module.  Denoting $(-)^{*}:=\Hom_{R}(-,R)\colon \mod R\to\mod R$, we say that $X\in\mod R$ is \emph{reflexive} if the natural map $X\to X^{**}$ is an isomorphism.

We denote $\refl R$ to be the category of reflexive $R$-modules, and we denote $\CM R$ to be the category of CM $R$-modules.  Throughout this section, we stress our convention that commutative rings are always assumed to be noetherian, and $R$ always denotes a commutative noetherian ring.

\subsection{$d$-CY and $d$-sCY algebras}\label{dCY}
When dealing with singularities, throughout this paper we use the language of $d$-sCY  algebras, introduced in \cite[\S3]{IR}.  Although technical, it provides the common language that links the commutative rings and the various noncommutative endomorphism rings that we will consider.

To do this, let $R$ be a commutative ring with $\dim R=d$ and let $\Lambda$ be a module-finite $R$-algebra (i.e.\ an $R$-algebra which is a finitely generated $R$-module).  For any $X\in\mod \Lambda$, denote by $E(X)$ the injective hull of $X$, and put $E:=\bigoplus_{\m\in\Max R}E(R/\m)$.  This gives rise to Matlis duality $D:=\Hom_R(-,E)$.  Note that if $R$ is over an algebraically closed field $\kk$, then $D$ coincides with $\Hom_{\kk}(-,\kk)$ on the category $\fl\Lambda$ by \cite[1.1, 1.2]{oo}. We let $\Dbf(\Lambda)$ denote all bounded complexes with finite length cohomology.
\begin{defin}
For $n\in\mathbb{Z}$, we call $\Lambda$ $n$-Calabi-Yau (=$n$-CY) if there is a functorial isomorphism 
\[
\Hom_{\D(\Mod \Lambda)}(X,Y[n])\cong D\Hom_{\D(\Mod \Lambda)}(Y,X)
\]
for all $X\in\Dbf(\Lambda)$ and $Y\in\Db(\mod\Lambda)$.  Similarly we call $\Lambda$ singular $n$-Calabi-Yau (=$n$-sCY)  if the above functorial isomorphism holds for all $X\in\Dbf(\Lambda)$ and $Y\in\Kb(\proj \Lambda)$.
\end{defin}

We will use the following characterizations.

\begin{prop}
\t{(1)} \cite[3.10]{IR} A commutative ring $R$ is $d$-sCY (respectively, $d$-CY) if and only if $R$ is Gorenstein (respectively, regular) and equi-codimensional, with $\dim R=d$.\\
\t{(2)} \cite[3.3(1)]{IR}  Let $R$ be a commutative $d$-sCY ring and $\Lambda$ a module-finite $R$-algebra which is a faithful $R$-module. Then $\Lambda$ is $d$-sCY if and only if $\Lambda\in\CM R$ and $\Hom_R(\Lambda,R)_\m\cong\Lambda_\m$ as $\Lambda_\m$-bimodules for all $\m\in\Max R$.
\end{prop}

Suppose that $R$ is a commutative $d$-sCY ring, and let $\Lambda$ be a module-finite $R$-algebra which is $d$-sCY.  We say that $X\in\mod\Lambda$ is \emph{maximal Cohen-Macaulay} (or simply, \emph{CM}) if $X\in\CM R$. We denote by $\CM\Lambda$ the category of CM $\Lambda$-modules. Then $X\in\mod\Lambda$ is CM if and only if $\Ext^{i}_{R}(X,R)=0$ for all $i\geq 1$, which is equivalent to $\Ext^{i}_{\Lambda}(X,\Lambda)=0$ for all $i\geq1$ by the proof of \cite[3.4(5)(i)]{IR}.  We denote $\uCM\Lambda$ to be the \emph{stable category}, where we factor out by those morphisms which factor through projective $\Lambda$-modules.  

On the other hand, when $\Lambda$ is a noetherian ring we denote
$$\Dsg(\Lambda):= \Db(\mod\Lambda)/\Kb(\proj\Lambda)$$
to be the \emph{singular derived} category of $\Lambda$. Since any $d$-sCY algebra satisfies $\id_\Lambda\Lambda=d$ and $\id_{\Lambda^{\op}}\Lambda=d$ by \cite[3.1(6)(2)]{IR}, we have the following equivalence by a standard theorem of Buchweitz \cite[4.4.1(2)]{Buch}.

\begin{thm}\label{Buchweitz}
Suppose that $R$ is a commutative $d$-sCY ring and $\Lambda$ is a module-finite $R$-algebra which is $d$-sCY, then there is a triangle equivalence
\[
\Dsg(\Lambda)\simeq \uCM\Lambda.
\] 
\end{thm}
By \ref{Buchweitz} we identify $\Dsg(\Lambda)$ and $\uCM\Lambda$ in the rest of this paper.  Now we recall two important notions.
\begin{defin}\cite{Auslanderisolated}\label{nonsingular}
Let $R$ be a commutative $d$-sCY ring and $\Lambda$ is a module-finite $R$-algebra which is $d$-sCY.\\
\t{(1)} We say that $\Lambda$ is \emph{non-singular} if $\gl\Lambda=d$.\\
\t{(2)} We say that $\Lambda$ has \emph{isolated singularities} if $\gl\Lambda_{\p}=\dim R_\p$ for all non-maximal primes $\p$ of $R$.
\end{defin}

\begin{prop}\label{Dbisolatedtest}
Suppose that $R$ is a commutative $d$-sCY ring and $\Lambda$ is a module-finite $R$-algebra which is $d$-sCY. Then\\
\t{(1)} $\Lambda$ is non-singular if and only if $\uCM\Lambda=0$.\\
\t{(2)} $\Lambda$ has isolated singularities if and only if all Hom-sets in $\uCM\Lambda$ are finite length $R$-modules. 
\end{prop}
\begin{proof}
These are well-known (e.g.\ \cite[3.3]{Y}, \cite{Auslanderisolated}), but for convenience of the reader we give the proof.\\
\t{(1)} See for example \cite[2.17(1)$\Leftrightarrow$(3)]{IW}.\\
\t{(2)} ($\Rightarrow$)  Since $\Hom_{\uCM \Lambda}(X,Y)\cong\Ext^1_{\Lambda}(\Omega^{-1} X,Y)$ and $\Omega^{-1} X\in\CM\Lambda$, we know that $\Supp_R\Hom_{\uCM \Lambda}(X,Y)$ consists of maximal ideals (e.g. \cite[2.6]{IW}). Thus the assertion follows.\\
($\Leftarrow$)  Let $\p$ be a prime ideal of $R$ with $\hgt\p<d$ and let $X\in\mod\Lambda_{\p}$.  Certainly we can find $Y\in\mod\Lambda$ with $Y_{\p}\cong X$ and so consider a projective resolution
\[
0\to K\to P_{d-1}\to \hdots\to P_{1}\to P_{0}\to Y\to 0
\]
in $\mod\Lambda$.  Since $\Lambda\in\CM R$, by localizing and using the depth lemma we see that $K\in\CM R$, i.e.\ $K\in\CM\Lambda$.  Consider
\begin{eqnarray}
0\to \Omega K\to P_{d}\to K\to 0\label{a},
\end{eqnarray}
then since $\Ext^{1}_{\Lambda}(K,\Omega K)\cong\Hom_{\uCM\Lambda}(\Omega K,\Omega K)$ has finite length, it is supported only on maximal ideals.  Hence (\ref{a}) splits under localization to $\p$ and so $K_{\p}$ is free.  Thus $\pd_{\Lambda_{\p}}(X)< d$.  By Auslander--Buchsbaum, $\gl\Lambda_{\p}=\dim R_\p$.
\end{proof}

\subsection{MM and CT modules}\label{subsection:MMandCT} 

Recall that $M\in \refl R$ gives an \emph{NCCR} $\Lambda:=\End_R(M)$ of $R$ if $\Lambda\in \CM R$ and $\gl\Lambda=d$ \cite{VdBNCCR}. In this case $\Lambda$ is $d$-sCY by \ref{ModisdCY}(2) below, and by definition $\Lambda$ is non-singular.

The following more general notions are quite natural from a
representation theoretic viewpoint.

\begin{defin}\label{MMdef}$\,$\cite{IW} Let $R$ be a $d$-sCY  ring.  Then\\
(1) $M\in\refl R$ is called a \emph{rigid module} if $\Ext^1_{R}(M,M)=0$.\\
(2) $M\in\refl R$ is called a \emph{modifying module} if $\End_{R}(M)\in\CM R$.\\
(3) We say $M\in\refl R$ is a \emph{maximal modifying (=MM) module} if it is modifying and further if $M\oplus Y$ is modifying for $Y\in\refl R$, then $Y\in\add M$.  Equivalently,
\[
\add M=\{X\in\refl R \mid \End_{R}(M\oplus X)\in\CM R\}.
\]
(4) We call $M\in\CM R$ a \emph{CT module} if
\[
\add M=\{X\in\CM R \mid \Hom_{R}(M,X)\in\CM R\}.
\]
\end{defin}

The following results will be used extensively.   As in \cite{IW}, if $X\in\mod R$, we denote $\fl X$ to be the largest finite length sub-$R$-module of $X$. 

\begin{lemma}\label{depthlemma}\cite[2.7, 5.12]{IW}  Suppose that $R$ is a 3-sCY  ring, and let $M\in\refl R$.\\
\t{(1)} If $M$ is modifying, then $\fl \Ext^1_R(M,M)=0$.  The converse holds if $M\in\CM R$.\\
\t{(2)} Assume that $R$ is an isolated singularity.  If $M$ is modifying, then it is rigid.  The converse holds if $M\in\CM R$.
\end{lemma}

\begin{lemma}\label{ModisdCY}
Suppose that $R$ is a d-sCY  normal domain, and $M\in\refl R$. Let $\Lambda=\End_R(M)$. Then\\
\t{(1)} \cite[2.4(3)]{IR} $\Lambda\cong\Hom_R(\Lambda,R)$ as $\Lambda$-bimodules.\\
\t{(2)} \cite[2.22(2)]{IW} $M$ is modifying if and only if $\Lambda$ is $d$-sCY . 
\end{lemma}

We remark that the property CT can be checked complete locally \cite[5.5]{IW}, that is $M$ is a CT $R$-module if and only if $\widehat{M}_\m$ is a CT $\widehat{R}_\m$-module for all $\m\in\Max R$.  In our context in \S\ref{Dsgresults}, this corresponds to the fact that a scheme $X$ is non-singular if and only if complete locally it is non-singular.  

In contrast, the property of being MM cannot be checked complete locally.   If $\widehat{M}_\m$ is an MM $\widehat{R}_\m$-module for all $\m\in\Max R$, then $M$ is an MM $R$-module.  However the converse is not true (see e.g.\ \ref{ex:MMnotcompletelocal} below). This corresponds to the difference between $\mathds{Q}$-factorial and complete locally $\mathds{Q}$-factorial singularities in \S\ref{Dsgresults}. 

We record here the following easy lemma, which shows the first link between modifying modules and factoriality.   

\begin{prop}\label{RisMMthenMMareFree}
Suppose $R$ is 3-sCY normal domain, with isolated singularities. Consider the statements\\
\t{(1)} $R$ is an MM module.\\
\t{(2)} Every modifying $R$-module is projective.\\
\t{(3)} $\uCM R$ is rigid-free.\\
\t{(4)} $R$ is locally factorial.\\
\t{(5)} $R$ is $\mathds{Q}$-factorial.\\
Then we have $\textnormal{(1)}\Leftrightarrow\textnormal{(2)}\Leftrightarrow\textnormal{(3)}\Rightarrow\textnormal{(4)}\Rightarrow\textnormal{(5)}$.
\end{prop}
\begin{proof}
Since $R$ is an isolated singularity, modifying modules are rigid by \ref{depthlemma}(2).\\
(1)$\Rightarrow$(2) Let $M\in\refl R$ be any modifying module.  Since $R$ is MM, there exists an exact sequence
\[
0\to F_{1}\to F_{0}\stackrel{f}{\to} M
\]
with $F_i\in \add R$ such that $f$ is a right ($\add R$)-approximation, by \cite[4.12]{IW}. Then clearly $f$ is surjective, and so we have $\pd_R(M)\leq 1$.  By \cite[4.10]{AG}, if $\pd_R(M) =1$ then $\Ext^1_R(M,M)\neq 0$, a contradiction.  Hence $M$ is projective.\\
(2)$\Rightarrow$(3) is clear, since $\Hom_{\uCM R}(M,M[1])\cong\Ext^1_R(M,M)$.\\
(3)$\Rightarrow$(1) Suppose that $\End_R(R\oplus X)\in\CM R$ for some $X\in\refl R$. Then $X=\Hom_R(R,X)$ is a CM $R$-module which is a rigid object in $\uCM R$. Hence $X$ is zero in $\uCM R$, i.e. $X\in\add R$.\\
(2)$\Rightarrow$(4) Since $R$ is a normal domain, all members of the class group have $R$ as endomorphism rings, and thus they are modifying $R$-modules.  Hence every member of the class group is projective, so the result follows.\\
(4)$\Rightarrow$(5) is clear.
\end{proof}

Below we will use the following result from commutative algebra. 
\begin{thm}$\,$\cite[3.1(1)]{DaoNCCR}\label{Dao}
Let $S$ be a regular local ring of dimension four containing a field $k$, and let $R=S/(f)$ be a hypersurface. If $R$ is a $\mathds{Q}$-factorial isolated singularity, then the free modules are the only modifying modules.
\end{thm}
\begin{proof}
This is the equi-characteristic version of \cite[3.1(1)]{DaoNCCR}, but since {\em loc.\ cit.\ }is in somewhat different language than used here, so we sketch Dao's proof.   Let $N$ be a modifying $R$-module.  Since $R$ is isolated, $N$ is rigid by \ref{depthlemma}(2).  On the other hand, by \cite{DaoRigid}, $N$ is a Tor-rigid $R$-module.  These two facts imply, via a result of Jothilingham \cite[Main Theorem]{Joth} (see also \cite[2.4]{DaoNCCR}), that $N$ is a free $R$-module.
\end{proof}

Now we have the following result, which strengthens \ref{RisMMthenMMareFree}.
\begin{thm}\label{UFDandrigid}
Let $R$ be a 3-sCY  normal domain over $\k$, which locally has only isolated hypersurface singularities.  Then the following are equivalent:\\
\t{(1)} $R$ is an MM module.\\
\t{(2)} Every modifying $R$-module is projective.\\
\t{(3)} $\uCM R$ is rigid-free.\\
\t{(4)} $R$ is locally factorial.\\
\t{(5)} $R$ is $\mathds{Q}$-factorial.
\end{thm}
\begin{proof}
By \ref{RisMMthenMMareFree} we only have to show (5)$\Rightarrow$(1).   Suppose that $M\in\refl R$ with $\End_R(R\oplus M)\in\CM R$.  Then for any $\m\in\Max R$ we have $\End_{R_\m}(R_\m\oplus M_\m)\in\CM R_\m$ with $R_\m$ satisfying the assumptions of \ref{Dao}, hence $M_\m$ is a free $R_\m$-module.  Thus $M_\m\in\add R_\m$ for all $\m\in\Max R$, so $M\in\add R$ \cite[2.26]{IW}.
\end{proof}

\begin{remark}
The corollary is false if we remove the isolated singularities assumption, since then (4) does not necessarily imply (1). An example is given by $R=\C{}[x,y,z,t]/(x^2+y^3+z^5)$.  Also, note that in the isolated singularity case it is unclear if the  hypersurface assumption is strictly necessary; indeed Dao conjectures that \ref{Dao} still holds in the case of complete intersections \cite[\S4]{DaoConj}.   We remark that removing the hypersurface assumption in \ref{UFDandrigid} would allow us to bypass the use of the classification of terminal singularities later.
\end{remark}

\begin{example}\label{ex:MMnotcompletelocal}
Consider the element $f=x^2+x^3+y^2-uv$ in the ring $\mathbb{C}[u,v,x,y]$, and set $R:=\mathbb{C}[u,v,x,y]/(f)$.  Then $R$ is an MM $R$-module, but $\widehat{R}_\m$ is not an MM $\widehat{R}_\m$-module, where $\m=(u,v,x,y)$. 
\end{example}
\begin{proof}
Since $R$ has an isolated singularity only at $\m:=(u,v,x,y)$
and $\sqrt{x+1}$ does not exist in $R_{\m}$, we have that
$R_\m$ is a factorial hypersurface singularity for all $\m\in \Max R$.  Hence if $M$ is a modifying $R$-module, then $M_\m$ is a free $R_\m$-module for any $\m\in \Max R$ by \ref{Dao}. Thus $M$ is projective.

For the last statement, since $\sqrt{x+1}$ exists (and is a unit) in the completion, $\widehat{R}_\m$ is isomorphic to $\mathbb{C}[[u,v,x,y]]/(uv-xy)$, for which $\widehat{R}_\m\oplus (u,x)$ is a modifying module. 
\end{proof}

Now we give a general categorical criterion for a given modifying module to be maximal, generalizing \ref{RisMMthenMMareFree}.

\begin{prop}\label{MMtestforuCM}
Suppose that $R$ is normal 3-sCY, let $M$ be a modifying $R$-module and set $\Lambda:=\End_{R}(M)$.  If $M$ is an MM $R$-module, then $\uCM\Lambda$ is rigid-free. Moreover the converse holds if $\Lambda$ has isolated singularities.
\end{prop}
\begin{proof}
By \cite[4.8(1)]{IW}, $M$ is an MM $R$-module if and only if there is no non-zero object $X\in\uCM\Lambda$ such that $\fl\Hom_{\uCM\Lambda}(X,X[1])=0$. This clearly implies $\uCM\Lambda$ is rigid-free. By \ref{Dbisolatedtest}(2), the converse is true if $\Lambda$ has isolated singularities.
\end{proof}

\begin{remark}
We conjectured in \cite{IW} that in dimension three MMAs always have isolated singularities (see \ref{mainQfact} for some evidence), so the key property from \ref{MMtestforuCM} is that the stable category of CM modules is rigid-free.  Note that $\Lambda:=\End_R(M)\in\CM R$ is an NCCR if and only if the stable category $\uCM\Lambda$ is zero.  Hence when passing from NCCRs to MMAs, the categories $\Dsg(\Lambda)\simeq\uCM \Lambda$ are no longer zero, but instead rigid-free. This motivates \S\ref{Dsgresults}.
\end{remark}

\subsection{Homologically finite complexes}
In this subsection we wish to consider the general setting of commutative rings, and so in particular the Krull dimension may at times be infinite.  The following is based on some Danish handwritten notes of Hans-Bj\o rn Foxby.
\begin{thm}\label{unbelievable}
Let $R$ be a commutative ring (not necessarily of finite Krull dimension), let $\Lambda$ be a module-finite $R$-algebra, and let $M\in\mod\Lambda$. If $\pd_{\Lambda_\p}(M_\p)<\infty$ for all $\p\in\Spec R$, then $\pd_\Lambda(M)<\infty$. 
\end{thm}

\begin{proof}
Take a projective resolution
\[
\hdots\to P_1\stackrel{d_1}{\to} P_0\stackrel{d_0}{\to}M\to 0,
\]
of $M$ with finitely generated projective $R$-modules $P_i$, and set $K_0:=\Im d_0=M$ and $K_i := \Im d_i=\Ker d_{i-1}$ for $i>0$.  Now
\[
\pd_{\Lambda}(M)< n  \iff \Ext^{n}_{\Lambda}(M,K_{n})= 0,
\]
and the right hand side is equivalent to the condition that $\Supp_R\Ext^{n}_{\Lambda}(M,K_{n})$ is empty. Localizing at $\p\in\Spec R$, we have
\begin{align*}
\pd_{\Lambda_\p}(M_\p)\geq n  &\iff \Ext^{n}_{\Lambda}(M,K_{n})_\p=\Ext^{n}_{\Lambda_\p}(M_\p,(K_{n})_\p)\neq 0\\
&\iff \p\in \Supp_R \Ext^{n}_{\Lambda}(M,K_{n}).
\end{align*} 
Since $\Ext^{n}_{\Lambda}(M,K_{n})$ is a finitely generated $R$-module, we have
\[
\Supp_R \Ext^{n}_{\Lambda}(M,K_{n})=V(\Ann_R(\Ext^{n}_{\Lambda}(M,K_{n}))). 
\]
For ease of notation set $I_n = \Ann_R( \Ext^{n}_{\Lambda}(M,K_{n}))$, then
\begin{eqnarray}
\pd_{\Lambda_\p}(M_\p)\geq n  \iff \p \in V(I_n).\label{u2}
\end{eqnarray}
In particular, we have a decreasing sequence
\[
V(I_0)\supseteq V(I_1)\supseteq V(I_2)\supseteq\hdots,
\]
so we have $V(I_t)=V(I_{t+1})=\hdots$ for some $t$ since the $\Spec R$ is noetherian by our assumption. If $\p\in V(I_t)=\bigcap_{n=t}^\infty V(I_n)$, then $\pd_{\Lambda_\p}(M_\p)=\infty$ by (\ref{u2}), a contradiction. Consequently $V(I_t)=\emptyset$, so we have $\Ext_\Lambda^{t}(M,K_{t})=0$ and $\pd_\Lambda(M)<t$.
\end{proof}

\begin{cor}\label{finitenumberofsimples} 
Let $R$ be a commutative ring and $\Lambda$ a module-finite $R$-algebra. Then \\
\t{(1)}  For every $\p\in\Spec R$, the algebra $\Lambda_{\p}$ has only finitely many simple modules.\\
\t{(2)}  Let $M\in\mod\Lambda$. Then $\pd_\Lambda(M)<\infty$ if and only if for all $X\in\mod\Lambda$ we have $\Ext_\Lambda^j(M,X)=0$ for $j\gg0$.
\end{cor}
\begin{proof}
(1) $\Lambda_{\p}$ is a module-finite $R_{\p}$-algebra.  There are only finitely many simple $\Lambda_{\p}$-modules since they are annihilated by $\p R_\p$ and hence they are modules over the finite dimensional $(R_\p/\p R_\p)$-algebra $\Lambda_\p/\p\Lambda_\p$, which has only finitely many simple modules.\\
(2) We only have to show `if' part. For any $\p\in\Spec R$, we know that the localization functor $(-)_\p\colon \mod\Lambda\to\mod\Lambda_\p$ is essentially surjective, and that the completion functor $\widehat{(-)}\colon \fl\Lambda_\p\to\fl\widehat{\Lambda}_\p$ is an equivalence since the completion does not change finite length modules. By (1) there exists $X\in\mod\Lambda$ such that $\widehat{X}_\p$ is the sum of all simple $\widehat{\Lambda}_\p$-modules. The assumptions imply that there exists $t\ge0$ such that $\Ext_\Lambda^j(M,X)=0$ for all $j>t$. Thus $\Ext^j_{\widehat{\Lambda}_\p}(\widehat{M}_\p,\widehat{X}_\p)=0$ for all $j>t$. This implies $\pd_{\widehat{\Lambda}_\p}(\widehat{M}_\p)\le t$ since $\widehat{M}_\p$ has a minimal projective resolution. Thus we have $\pd_{\Lambda_\p}(M_\p)\le t$ by applying \cite[2.26, 2$\Leftrightarrow$4]{IW} to the $t$-th syzygy of $M$. By \ref{unbelievable}, it follows that $\pd_\Lambda(M)<\infty$.
\end{proof}

Let $\c{T}$ be a triangulated category with a suspension functor $[1]$.  Recall that $x\in \c{T}$ is called {\em homologically finite} if for all $y\in \c{T}$, $\Hom_{\c{T}}(x,y[i])=0$ for all but finitely many $i$.  The following is well-known under more restrictive hypothesis:

\begin{prop}\label{homfinitecomplexesOK}
Let $R$ be a commutative ring and $\Lambda$ be a module-finite $R$-algebra.  Then the homologically finite complexes in $\Db(\mod\Lambda)$ are precisely $\Kb(\proj\Lambda)$.
\end{prop}
\begin{proof}
Since every object in $\Kb(\proj \Lambda)$ is homologically finite, we just need to show the converse.  If $X$ is homologically finite, replace $X$ by its projective resolution $(P,d)$.  The truncation of $P$, namely $\tau_{\geq i}P$, is quasi-isomorphic to $X$ for small enough $i$.  Fix such an $i$.  Thus, if we can show that $\Ker(d_i)$ has finite projective dimension, certainly it follows that $X$ belongs to $\Kb(\proj \Lambda)$.

Now by the usual short exact sequence given by truncation, there is a triangle 
\[
Q \to X \to \Ker(d_i)[i] \to Q[1]
\]
where $Q\in\Kb(\proj\Lambda)$.  Since the first two are homologically finite, so is $\Ker(d_i)$.  In particular, for all $Y\in\mod \Lambda$, we have that $\Ext_\Lambda^j(\Ker(d_i),Y)=0$ for large enough $j$.  By \ref{finitenumberofsimples}(2), $\Ker(d_i)$ has finite projective dimension, as required.  
\end{proof}

\section{Singular derived categories and $\mathds{Q}$-factorial terminalizations}\label{Dsgresults}

\subsection{Singular derived categories of Gorenstein schemes with isolated singularities}\label{Dsgprelim}
Recall that a scheme $X$ is called {\em Gorenstein} if $\c{O}_{X,x}$ is a Gorenstein local ring for all closed points $x\in X$.  We say that a scheme $X$ satisfies \emph{(ELF)} if $X$ is separated, noetherian, of finite Krull dimension, such that $\coh X$ has enough locally free sheaves (i.e.\ for all $\c{F}\in\coh X$, there exists a locally free sheaf $\c{E}$ and a surjection $\c{E}\twoheadrightarrow\c{F}$). This is automatic in many situations, for example if $X$ is quasi--projective over a commutative noetherian ring \cite[2.1.2(c), 2.1.3]{TT}.

\begin{defin}
We say that $\F\in\Qcoh X$ is \emph{locally free} if $X$ can be covered with open sets $U$ for which each $\F|_U$  is a (not necessarily finitely generated) free $\O_X|_U$-module.
\end{defin}
If $\F\in\Qcoh X$ is locally free, then $\F_x$ is a free
$\O_{X,x}$-module for all closed points $x\in X$. The converse
is true if $\F\in\coh X$.

As in \cite[\S1]{Orlov} we consider
\begin{itemize}
\item the full triangulated subcategory $\Perf(X)$ of $\Db(\coh X)$ consisting of objects which are isomorphic to bounded complexes of finitely generated locally free sheaves in $\Db(\coh(X))$. 
\item the full triangulated subcategory $\Lfr(X)$ of $\Db(\Qcoh X)$ consisting of objects which are isomorphic to bounded complexes of locally free sheaves in $\Db(\Qcoh(X))$. 
\end{itemize}
Following Orlov \cite{Orlov,OrlovCompletion}, we denote
$$\Dsg(X):=\Db(\coh X)/\Perf(X)\ \ \ \mbox{ and }\ \ \ \DsgQ(X):=\Db(\Qcoh X)/\Lfr(X).$$ 
Then the natural functor $\Dsg(X)\to\DsgQ(X)$ is fully faithful \cite[1.13]{Orlov}, and we regard $\Dsg(X)$ as a full subcategory of $\DsgQ(X)$.

The aim of this subsection is to prove the following result, which plays a crucial role in this paper.  Note that part (2) also follows from work of Orlov \cite{OrlovCompletion}, and is well-known to experts \cite{BK}.
\begin{thm}\label{Dsgsplits}
Suppose that $X$ is a Gorenstein scheme of dimension $d$ satisfying (ELF), with isolated singularities $\{x_1,\hdots,x_n\}$. \\
\t{(1)} There is a triangle equivalence $\Dsg(X)\to\bigoplus_{i=1}^{n}\Dsg(\c{O}_{X,x_i})=\bigoplus_{i=1}^{n}\uCM \c{O}_{X,x_i}$ up to summands,
given by $\F\mapsto (\F_{x_1},\hdots,\F_{x_n})$.\\
\t{(2)} Taking the idempotent completion gives a triangle equivalence
$\overline{\Dsg(X)}\simeq\bigoplus_{i=1}^{n}\uCM\widehat{\c{O}}_{X,x_i}$.
\end{thm}

Since a key role in our proof of Theorem \ref{Dsgsplits} is played by the category $\DsgQ(X)$, we need to deal with infinitely generated modules. Let us start with recalling results by Gruson--Raynaud \cite{GR} and Jensen \cite{Jensen}.

\begin{prop}\label{prepare GR and J}
Let $R$ be a noetherian ring of finite Krull dimension and $M\in\Mod R$. Then\\
\t{(1)} $\FPD(R):=\sup\{ \pd(M)\mid M\in\Mod{R}\text{ and } \pd_R(M)<\infty \}=\dim R$.\\
\t{(2)} $\fd_R(M)<\infty$ implies $\pd_R(M)<\infty$.
\end{prop}
\begin{proof}
(1) is due to Bass and Gruson--Raynaud \cite[II.3.2.6]{GR} (see also \cite[3.2]{Foxby}).\\
(2) is due to Jensen \cite[Prop.6]{Jensen} (see also \cite[3.3]{Foxby}).
\end{proof}

We have the following immediate consequence.

\begin{lemma}\label{finite_proj_dim}
Let $R$ be a noetherian ring of finite Krull dimension. Then for any $M\in\Mod R$, $\pd_{R_{\m}}(M_{\m})< \infty\text{ for all }\m\in\Max R \iff \pd_R(M)\leq \dim R$.
\end{lemma}
\begin{proof}
($\Leftarrow$) is trivial\\
($\Rightarrow$) The hypothesis implies that $\pd_{R_{\p}}(M_{\p})<\infty$ for all $\p\in\Spec R$. Thus by \ref{prepare GR and J}(1), we have $\pd_{R_{\p}}(M_{\p})\leq \dim R_{\p}\leq \dim R$ for all $\p\in\Spec R$.  In particular $\fd_{R_\p}(M_\p)\leq \dim R$ for all $\p\in\Spec R$.  Now $X\in \Mod R$ is zero if $X_\p=0$ for all $\p\in \Spec R$. This implies $\fd_R(M)\le\dim R$ since Tor groups localize also for infinitely generated modules. By \ref{prepare GR and J}(2), $M$ has finite projective dimension. Again by \ref{prepare GR and J}(1), we have $\pd_R(M)\le\dim R$.
\end{proof}

The following is easy:
\begin{lemma}\label{KbProj=KbFree} 
Let $R$ be a commutative noetherian ring of finite Krull dimension, then $\Kb(\Free R)=\Lfr(\Spec R)=\Kb(\Proj R)$.
\end{lemma}
\begin{proof}
We have natural inclusions $\Kb(\Free R)\subseteq \Lfr(\Spec R)\subseteq \Kb(\Proj R)$ since every free module is clearly locally free, and further every locally free $R$-module is locally projective, thus projective by \cite[II.3.1.4(3)]{GR}. We only have to show $\Kb(\Proj R)=\Kb(\Free R)$. The proof is similar to \cite[2.2]{Rickard} --- by the Eilenberg swindle, if $P\in\Proj R$, there exists $F\in\Free R$ such that $P\oplus F$ is free.  Hence from every object of $\Kb(\Proj R)$ we can get to an object of $\Kb(\Free R)$ by taking the direct sum with complexes of the form $0\to F\to F\to 0$, which are zero objects in the homotopy category.
\end{proof}

\begin{cor}\label{fpdgivesLfr}
Suppose that $X$ satisfies (ELF). If $\F\in\Qcoh X$ satisfies $\pd_{\c{O}_{X,x}}(\F_x)<\infty$ for all closed points $x\in X$, then $\F\in\Lfr(X)$.
\end{cor}
\begin{proof}
Since $X$ satisfies (ELF), membership of $\Lfr(X)$ can be checked locally (see \cite[1.7, proof of 1.14]{Orlov}). Thus we only have to show $\F|_U\in\Lfr(U)$ for each affine open subset $U=\Spec R$ of $X$. Since $\F_{\m}$ is an $R_\m$-module with finite projective dimension for all $\m\in\Max R$, we have $\F|_U\in\Kb(\Proj R)$ by \ref{finite_proj_dim}. By \ref{KbProj=KbFree}, we have $\F|_U\in\Lfr(U)$, so the assertion follows.
\end{proof}

We will also require the following well-known lemma, due to Orlov.  If $s$ is an object in a triangulated category $\c{T}$, we denote $\thick_{\c{T}}(s)$ to be the smallest triangulated subcategory of $\c{T}$ which is closed under direct summands and isomorphisms and contains $s$.
\begin{lemma}\t{(Orlov)}\label{isolatedGEN}
Suppose that $X$ satisfies (ELF), and has only isolated singularities $\{ x_1,\hdots,x_n\}$.   Denote the corresponding skyscraper sheaves by $k_1,\hdots,k_n$.  Then $\Dsg(X)=\thick_{\Dsg(X)}(\bigoplus_{i=1}^nk_i)$.
\end{lemma}
\begin{proof}
By assumption the singular locus consists of a finite number of closed points, hence is closed.  It follows that $\Dsg(X)=\thick_{\Dsg(X)}(\coh_{\{x_1,\hdots,x_n \}}X)$ by \cite[1.2]{ChenOrlov}.  It is clear that every coherent sheaf supported in $\{x_1,\hdots,x_n \}$ belongs to $\thick_{\Dsg(X)}(\bigoplus_{i=1}^nk_i)$.
\end{proof}

Now we are ready to prove \ref{Dsgsplits}.

For each $x_i$, consider a morphism $f_i:=(\Spec \c{O}_{X,x_i}\stackrel{g_i}{\to}\Spec R_i\stackrel{h_i}{\to} X)$ where $\Spec R_i$ is some affine open subset of $X$ containing $x_i$. Let $S:=\c{O}_{X,x_1}\oplus\hdots\oplus \c{O}_{X,x_n}$ and $Y:=\Spec S$. Then the collection of the $f_i$ induce a morphism
\[Y=\textstyle\coprod_{i=1}^{n}\Spec \c{O}_{X,x_i}\stackrel{f}{\to} X.\]
The functor $g_{i\,*}$ is just extension of scalars corresponding to a localization $R_i\to {R_i}_\m$, so $g_{i\,*}$ is exact and preserves quasi-coherence.  Further $h_i$ is an affine morphism (since $X$ is separated), hence $h_{i\,*}$ also is exact and preserves quasi-coherence.  Hence each $f_{i\, *}$ is exact and preserves quasi-coherence.  
Thus we have an adjoint pair
\[
\begin{tikzpicture}
\node (A) at (-0.75,0) {$\Qcoh X$};
\node (B) at (3.25,0) {$\Qcoh Y=\bigoplus\limits_{i=1}^{n}\Qcoh \c{O}_{X,x_i}$};
\draw[->] (0,0.1) -- node [above] {$\scriptstyle f^*$} (1,0.1);
\draw [<-] (0,-0.1) -- node [below] {$\scriptstyle f_*$} (1,-0.1);
\end{tikzpicture}
\]
which are explicitly given by $f^*\c{F}=(\c{F}_{x_1},\hdots,\c{F}_{x_n})$ and $f_*(\c{G}_1,\hdots,\c{G}_n)=\bigoplus_{i=1}^{n} f_{i\, *}\c{G}_i$.  These functors are exact, so induce an adjoint pair
\[
\begin{tikzpicture}
\node (A) at (-1.25,0) {$\Db(\Qcoh X)$};
\node (B) at (4,0) {$\Db(\Qcoh Y)=\bigoplus\limits_{i=1}^{n}\Db(\Qcoh \c{O}_{X,x_i})$};
\draw[->] (0,0.1) -- node [above] {$\scriptstyle f^*$} (1,0.1);
\draw [<-] (0,-0.1) -- node [below] {$\scriptstyle f_*$} (1,-0.1);
\end{tikzpicture}
\]
in the obvious way. 

We will show $f^*(\Lfr(X))\subseteq\Lfr(Y)$ and $f_*(\Lfr(Y))\subseteq\Lfr(X)$, since this then induces an adjoint pair
\[
\begin{tikzpicture}
\node (A) at (-0.75,0) {$\DsgQ(X)$};
\node (B) at (3.25,0) {$\DsgQ(Y)=\bigoplus_{i=1}^{n}\DsgQ(\c{O}_{X,x_i})$.};
\draw[->] (0,0.1) -- node [above] {$\scriptstyle f^*$} (1,0.1);
\draw [<-] (0,-0.1) -- node [below] {$\scriptstyle f_*$} (1,-0.1);
\end{tikzpicture}
\]
by \cite[1.2]{Orlov}. It is clear that $f^*$ takes locally free sheaves to projective modules, and so $f^*(\Lfr(X))\subseteq\Lfr(Y)$. On the other hand, each $f_{i\, *}$ is an affine, flat morphism so $f_{i\, *}\c{O}_{X,x_i}$ is a flat $\c{O}_X$-module.  Further, each $f_{i\, *}$ preserves sums, so takes projective $\c{O}_{X,x_i}$-modules to flat $\c{O}_{X}$-modules.  Consequently, it follows that $f_*$ takes a projective $S$-module $P$ to a sheaf $f_*(P)$ for which $f_*(P)_x$ is a flat $\c{O}_{X,x}$-module for all closed points $x\in X$. By \ref{prepare GR and J}(2) we have $\pd_{\c{O}_{X,x}}(f_*(P)_x)<\infty$ for all closed points $x\in X$, hence by \ref{fpdgivesLfr} $f_*(P)\in\Lfr(X)$ holds. Thus we have $f_*(\Lfr(Y))\subseteq\Lfr(X)$.

Now denote the skyscraper sheaves in $X$ corresponding to the singular points $x_1,\hdots,x_n$ by $k_{1},\hdots, k_{n}$.  By \ref{isolatedGEN} we know that $\Dsg(X)=\thick_{\Dsg(X)}(\bigoplus_{i=1}^n k_{i})$. This implies $\Dsg(X)\subseteq\thick_{\DsgQ(X)}(\bigoplus_{i=1}^n k_{i})$.

Let $\alpha\colon 1\to f_* \circ f^*$ be the unit. Since $f_* \circ f^*(k_i)=f_* (0,\hdots,0,k_i,0,\hdots,0)=f_{i\, *} k_i=k_i$ for all $1\leq i\leq n$, we have that $\alpha\colon 1\to f_* \circ f^*$ is an isomorphism on $\thick_{\DsgQ(X)}(\bigoplus_{i=1}^n k_{i})$. In particular $f^*$ is fully faithful on $\thick_{\DsgQ(X)}(\bigoplus_{i=1}^n k_{i})$. Since $f^*$ clearly takes $\Dsg(X)$ to $\Dsg(Y)$, and $\Dsg(X)\subseteq\thick_{\DsgQ(X)}(\bigoplus_{i=1}^n k_{i})$, we deduce that 
\[
\Dsg(X)\stackrel{f^*}{\longrightarrow}\Dsg(Y)
\]
is fully faithful. 
On the other hand, since $\Dsg(Y)=\thick_{\Dsg(Y)}(f^*\bigoplus_{i=1}^n k_{i})$ holds by \ref{isolatedGEN}, we have $\Dsg(Y)=\thick_{\Dsg(Y)}(f^*(\Dsg(X)))$. This immediately implies that $f^*\colon \Dsg(X)\to\Dsg(Y)$ is an equivalence up to direct summands (e.g. \cite[2.1.39]{Neeman}).

The first statement (1) in the theorem now follows from the well-known equivalence $\Dsg(\c{O}_{X,x_i})\simeq \uCM\c{O}_{X,x_i}$ due to Buchweitz \cite[4.4.1]{Buch}.  The second statement (2) follows from (1), since the idempotent completion of $\uCM\c{O}_{X,x_i}$ is $\uCM\c{\widehat{O}}_{X,x_i}$ \cite[A.1]{KMV}.
\qed

\medskip
The first corollary of \ref{Dsgsplits} is the following alternative proof of \cite[1.24]{Orlov}, which we will use later:
\begin{cor}\label{isolatedsingsfinitelength}
Suppose that $X$ is a $d$-dimensional Gorenstein scheme over $\k$, satisfying (ELF), with isolated singularities.  Then all Hom-sets in $\Dsg(X)$ are finite dimensional $k$-vector spaces.
\end{cor}
\begin{proof}
By \ref{Dsgsplits}(1), it is enough to show that all Hom-sets in $\uCM\c{O}_{X,x_i}$ are finite dimensional $k$-vector spaces. Since each $\c{O}_{X,x_i}$ is isolated, all Hom-sets in $\uCM\c{O}_{X,x_i}$ are finite length $\c{O}_{X,x_i}$-modules by \ref{Dbisolatedtest}(2). But each $\c{O}_{X,x_i}$ is a localization of a finitely generated $\k$-algebra, so by the Nullstellensatz its residue field is a finite extension of $\k$. Thus the assertion follows.
\end{proof}

\subsection{Rigid-freeness and $\mathds{Q}$-factoriality}
In this section we give our main result which characterizes the $\mathds{Q}$-factorial property in terms of the singular derived category, and then relate this to MMAs. 

\begin{defin}
Suppose that $X$ is a Gorenstein scheme.  We say that $\F\in\coh X$ is a \emph{Cohen--Macaulay} (=CM) sheaf if $\F_x$ is a CM $\c{O}_{X,x}$-module for all closed points $x\in X$.
\end{defin}
Under the assumption that $X$ is Gorenstein, $\F\in\coh X$ is a CM sheaf if and only if  
\[
\c{E}xt^i(\F,\c{O}_X)_x(=\Ext^i_{\c{O}_{X,x}}(\F_x,\c{O}_{X,x}))=0
\]
for all closed points $x\in X$ and all $i>0$. Thus $\F\in\coh X$ is a CM sheaf if and only if $\c{E}xt^i(\F,\c{O}_X)=0$ for all $i>0$.

The following is a generalization of \ref{RisMMthenMMareFree}(3)$\Rightarrow$(4).
\begin{prop}\label{locallyfactorial}
Suppose that $X$ is a Gorenstein normal scheme satisfying (ELF), of dimension three which has only isolated singularities $\{ x_1,\hdots,x_n\}$.  If $\Dsg(X)$ is rigid-free, then $X$ is locally factorial.
\end{prop}
\begin{proof}
We need to show that any reflexive sheaf $\c{G}$ on $X$ of rank one is locally free. For all closed points $x\in X$, we have a reflexive $\c{O}_{X,x}$-module $\c{G}_{x}$ of rank one. Since $\c{O}_{X,x}$ is a normal domain, we have $\End_{\c{O}_{X,x}}(\c{G}_x)=\c{O}_{X,x}$. Since $\c{O}_{X,x}$ is an isolated singularity, we have $\Ext^1_{\c{O}_{X,x}}(\c{G}_x,\c{G}_x)=0$ by \ref{depthlemma}(2). 

Since $X$ satisfies (ELF), there exists an exact sequence $0\to\c{K}\to\c{V}\to\c{G}\to 0$ with a locally free sheaf $\c{V}$. Localizing at $x$, we have an exact sequence
\begin{eqnarray}
0\to\c{K}_x\to\c{V}_x\to\c{G}_x\to 0\label{localses}
\end{eqnarray}
of $\c{O}_{X,x}$-modules with a projective $\c{O}_{X,x}$-module $\c{V}_x$. Since $\c{G}_x\in\refl\c{O}_{X,x}$, we have $\c{K}_x\in\CM\c{O}_{X,x}$. Since $\Ext^1_{\c{O}_{X,x}}(\c{G}_x,\c{G}_x)=0$, we have $\End_{\c{O}_{X,x}}(\c{K}_x)\in\CM\c{O}_{X,x}$ by applying \cite[4.10]{IW} to (\ref{localses}). Thus $\Ext^1_{\c{O}_{X,x}}(\c{K}_{x},\c{K}_{x})=0$ holds, again by \ref{depthlemma}(2).

Since $\c{K}_x\in\CM\c{O}_{X,x}$, we have by \ref{Dsgsplits}(1)
\[\Hom_{\Dsg(X)}(\c{K},\c{K}[1])\cong\bigoplus_{i=1}^{n}\Ext^1_{\c{O}_{X,x_i}}(\c{K}_{x_i},\c{K}_{x_i})=0.\]
Since $\Dsg(X)$ is rigid-free by our assumption, we have $\c{K}\in\Perf(X)$. Thus $\c{K}_x$ is a CM $\c{O}_{X,x}$-module of finite projective dimension. By the Auslander--Buchsbaum equality, $\c{K}_x$ has to be a projective $\c{O}_{X,x}$-module. By (\ref{localses}), we have $\pd_{\c{O}_{X,x}}\c{G}_x\leq 1$ for all $x\in X$. If $\pd_{\c{O}_{X,x}}\c{G}_x=1$, then $\Ext^1_{\c{O}_{X,x}}(\c{G}_x,\c{G}_x)\neq 0$ by \cite[4.10]{AG}, a contradiction. Thus $\c{G}_x$ is a projective $\c{O}_{X,x}$-module for all $x\in X$. Hence $\c{G}$ is locally free.
\end{proof}

Now we are ready to state our main result, which gives a relationship between factoriality of schemes and rigid-freeness of their singular derived categories.

\begin{thm}\label{Section3_main} Suppose that $X$ is a normal 3-dimensional Gorenstein scheme over $\k$, satisfying (ELF), with isolated singularities $\{x_1,\hdots,x_n\}$. \\
\t{(1)} If $\c{O}_{X,x_i}$ are hypersurfaces for all $1\leq i\leq n$,  then the following are equivalent:
\begin{enumerate}
\item[(a)] $X$ is locally factorial.
\item[(b)] $X$ is $\mathds{Q}$-factorial.
\item[(c)] $\Dsg(X)$ is rigid-free.
\item[(d)] $\uCM\O_{X,x}$ is rigid-free for all closed points $x\in X$.
\end{enumerate}
\t{(2)} If $\widehat{\c{O}}_{X,x_i}$ are hypersurfaces for all $1\leq i\leq n$,  then the following are equivalent: 
\begin{enumerate}
\item[(a)] $X$ is complete locally factorial.
\item[(b)] $X$ is complete locally $\mathds{Q}$-factorial.
\item[(c)] $\overline{\Dsg(X)}$ is rigid-free.
\item[(d)] $\uCM\widehat{\O}_{X,x}$ is rigid-free for all closed points $x\in X$.
\end{enumerate}
\end{thm}
\begin{proof}
(1) (a)$\Rightarrow$(b) is clear.  (b)$\Rightarrow$(d) follows from \ref{UFDandrigid}(5)$\Rightarrow$(3). (d)$\Rightarrow$(c) follows from \ref{Dsgsplits}(1). 
(c)$\Rightarrow$(a) follows from \ref{locallyfactorial}.\\
(2) The proof is identical to the proof in (1).
\end{proof}

\section{MMAs and NCCRs via tilting}\label{sect4}

The aim of this section is to give explicit information on rings that are derived equivalent to certain varieties.  Our main results (\ref{Endincorrectform}, \ref{CMiffcrepant}, \ref{NCCRdimd} and \ref{mainQfact}) show a strong relationship between the crepancy of birational morphisms and the property that the rings are Cohen-Macaualy.  In particular, we show that all algebras derived equivalent to $\mathds{Q}$-factorial terminalizations in dimension three are MMAs.  

Throughout this section, we continue to use the word variety to mean a normal integral scheme, of finite type over a field $k$, satisfying (ELF).  Some of the results below remain true with weaker assumptions.  If $X$ is a variety, we denote $\cd_X:=\gsh k$ to be the dualizing complex of $X$ \cite[Ch V \S10]{RD}, where $g\colon X\to\Spec k$ is the structure morphism.  If further $X$ is CM of dimension $d$, we will always use $\omega_X$ to denote $\cd_X[-d]$, and refer to it as the {\em geometric canonical}.  Thus, although canonical modules and dualizing complexes are not unique, throughout this section $\omega_X$ has a fixed meaning.  Note that crepancy is defined with respect to this canonical, i.e.\ $f\colon Y\to X$ is called {\em crepant} if $f^*\omega_X=\omega_Y$.

\subsection{The endomorphism ring of a tilting complex}

We call $\c{V}\in\Perf(Y)$ a \emph{pretilting complex} if $\Hom_{\D(Y)}(\c{V},\c{V}[i])=0$ for all $i\neq 0$, and a \emph{tilting complex} if it is a pretilting complex satisfying $\thick_{\D(Y)}(\c{V})=\Perf(Y)$.   If $Y$ is derived equivalent to $\Lambda$, there exists a tilting complex $\c{V}$ such that $\Lambda\cong\End_{\D(Y)}(\c{V})$.

The following two results are well-known.

\begin{lemma}\label{DbinducesPerf}
Suppose that $X$ is a scheme over a field $k$ satisfying (ELF).  Suppose that $X$ is derived equivalent to some ring $\Lambda$, then\\
\t{(1)} $\Perf(X)\simeq\Kb(\proj\Lambda)$.\\
\t{(2)} $\Dsg(X)\simeq\Dsg(\Lambda)$.
\end{lemma}
\begin{proof}
(1) Clearly the derived equivalence $\Db(\coh X)\simeq\Db(\mod\Lambda)$ induces an equivalence between full subcategories of homologically finite complexes. We observed in \ref{homfinitecomplexesOK} that $\Kb(\proj\Lambda)$ are precisely the homologically finite complexes in $\Db(\mod\Lambda)$. On the other hand, $\Perf(X)$ are precisely the homologically finite complexes in $\Db(\coh X)$ by \cite[1.11]{Orlov2}. Thus the assertion follows.\\
(2) follows immediately from (1).
\end{proof}

\begin{cor}\label{Dbisolated}
Suppose that $X$ is a scheme over a field $k$ satisfying (ELF).  Suppose that $X$ is derived equivalent to some ring $\Lambda$, which is a $d$-sCY algebra over a normal $d$-sCY $\k$-algebra. Then\\
\t{(1)} $X$ is smooth if and only if $\Lambda$ is non-singular (in the sense of \ref{nonsingular}).\\
\t{(2)} If $X$ has only isolated singularities, then $\Lambda$ has isolated singularities.
\end{cor}
\begin{proof}
The derived equivalence induces $\Dsg(X)\simeq\Dsg(\Lambda)$ by \ref{DbinducesPerf}.  It is very well-known that $X$ is smooth if and only if $\Dsg(X)=0$, so (1) follows from \ref{Dbisolatedtest}(1).  (2) is immediate from \ref{isolatedsingsfinitelength} and \ref{Dbisolatedtest}(2).
\end{proof}

In the rest of this subsection we suppose that we have a projective birational morphism $f\colon Y\to \Spec R$ and that $Y$ is derived equivalent to some ring $\Lambda$, and we investigate when $\Lambda$ has the form $\End_R(M)$ for some $M\in\refl R$.  We need the following well-known lemma.

\begin{lemma}\label{FBCtilt}
Let $Y\to\Spec R$ be a projective birational morphism between $d$-dimensional normal integral schemes. Let $\p\in\Spec R$ and consider the following pullback diagram:
\[
{\SelectTips{cm}{10}
\xy
(20,0)*+{Y^\prime}="2",(45,0)*+{Y}="3",
(20,-12.5)*+{\Spec R_\p}="b2",(45,-12.5)*+{\Spec R}="b3"
\ar^{i}"2";"3"
\ar^{j}"b2";"b3"
\ar^{g}"2";"b2"
\ar^{f}"3";"b3"
\endxy}
\]
\t{(1)} If $\p$ is a height one prime, then $g$ is an isomorphism.\\
\t{(2)} If $\c{V}$ is a pretilting complex of $Y$ with $\Lambda:=\End_{\D(Y)}(\c{V})$, then $\Lambda$ is a module-finite $R$-algebra  and $i^*\c{V}$ is a pretilting complex of $Y^\prime$ with $\Lambda_\p\cong\End_{Y^\prime}(i^*\c{V})$.
\end{lemma}
\begin{proof}
(1) Since $Y$ and $\Spec R$ are integral schemes and $f$ is a projective birational morphism, it is well-known that $f$ is the blowup of some ideal $I$ of $R$ \cite[8.1.24]{Liu}.  But blowups are preserved under flat base change, and in fact $g$ is the blowup at the ideal $I_\p$ \cite[8.1.14]{Liu}.  Now since $R$ is normal $R_\p$ is a discrete valuation ring, so every ideal is principal.  Thus $g$ is the blowup of a Cartier divisor, and so is an isomorphism by the universal property of blowing up.\\
(2)  Since $\coh Y$ is an $R$-linear category, so is its derived category,  and hence $\Lambda$ has the structure of an $R$-algebra. Since $\c{V}$ is a pretilting complex, we have
\begin{eqnarray}
\Lambda=\RHom_Y(\c{V},\c{V})=\Rf\RsHom_Y(\c{V},\c{V}).\label{last eq}
\end{eqnarray}
Now since $\c{V}$ is perfect, $\RsHom_Y(\c{V},\c{V})$ is a bounded complex of coherent sheaves.  Further, since $f$ is proper, $f_*$ preserves coherence and so $\Rf\colon\Db(\coh Y)\to \Db(\mod R)$. Therefore $\Rf\RsHom_Y(\c{V},\c{V})=\Lambda$ is a finitely generated $R$-module.

Lastly, applying $j^*$ to both sides of \eqref{last eq} gives
\begin{eqnarray*}
\Lambda_\p&=&j^*\Rf\RsHom_Y(\c{V},\c{V})=\Rg i^*\RsHom_Y(\c{V},\c{V})=\Rg \RsHom_{Y^\prime}(i^*\c{V},i^*\c{V})\\
&=&\RHom_{Y^\prime}(i^*\c{V},i^*\c{V}),
\end{eqnarray*}
where the second equality is flat base change, and the third holds since $i$ is flat and $\c{V}$ is coherent \cite[II.5.8]{RD}. Thus we have the assertion.
\end{proof}

We will also need the following result of Auslander--Goldman.  A proof can be found in \cite[2.11]{IW}.

\begin{prop}\cite{AG}\label{Auslander-Goldman}
Let $R$ be a normal domain, and let $\Lambda$ be a module-finite $R$-algebra. Then the following conditions are equivalent:\\
\t{(1)} There exists $M\in\refl R$ such that $\Lambda\cong\End_R(M)$ as $R$-algebras.\\
\t{(2)} $\Lambda\in\refl R$ and further $\Lambda_\p$ is Morita equivalent to $R_\p$ for all $\p\in\Spec R$ with $\hgt\p=1$.
\end{prop}
This leads to the following, which is an analogue of \cite[4.6(1)]{IW}.

\begin{thm}\label{Endincorrectform}
Let $Y\to\Spec R$ be a projective birational morphism between $d$-dimensional normal integral schemes, and suppose that $Y$ is derived equivalent to some ring $\Lambda$.  Then $\Lambda$ is a module-finite $R$-algebra.  If moreover $\Lambda\in\refl R$, then $\Lambda\cong\End_R(M)$ for some $M\in\refl R$.
\end{thm}
\begin{proof}
Since $\Lambda\cong\End_{\D(Y)}(\c{V})$ for some tilting complex $\c{V}$, it is a module-finite $R$-algebra by \ref{FBCtilt}(2).  Now consider a height one prime $\p\in\Spec R$, then by base change we have a pullback diagram
\[
{\SelectTips{cm}{10}
\xy
(20,0)*+{Y^\prime}="2",(45,0)*+{Y}="3",
(20,-12.5)*+{\Spec R_\p}="b2",(45,-12.5)*+{\Spec R}="b3"
\ar^{i}"2";"3"
\ar^{j}"b2";"b3"
\ar^{g}"2";"b2"
\ar^{f}"3";"b3"
\endxy}
\]
By \ref{FBCtilt}(1) $g$ is an isomorphism, and by \ref{FBCtilt}(2) $\Lambda_\p=\RHom_{Y^\prime}(i^*\c{V},i^*\c{V})$. Note that $i^*\c{V}\neq 0$ since $\Lambda\in\refl R$, so since $R$ is normal necessarily $\Lambda$ is supported everywhere.  Since $R_\p$ is a local ring, the only perfect complexes $x$ with $\Hom_{\D(R_\p)}(x,x[i])=0$ for all $i>0$ are shifts of projective modules \cite[2.12]{RZ}.   Thus $g_*i^*\c{V}\cong R_\p^a[b]$ for some $a\in\mathbb{N}$ and $b\in\mathbb{Z}$.  Hence $\Lambda_\p\cong\End_{R_\p}(R_\p^a)$, which is Morita equivalent to $R_\p$.  This holds for all height one primes, so by \ref{Auslander-Goldman} the assertion follows.
\end{proof}

\subsection{Cohen--Macaulayness and crepancy}

Let $f\colon Y\to\Spec R$ be a projective birational morphism such that $Y$ is derived equivalent to some ring $\Lambda$. In this section we show that $f$ is crepant if and only if $\Lambda\in\CM R$. To do this requires the following version of Grothendieck duality.  
\begin{thm}\label{SDT}
Suppose that $f\colon Y\to \Spec R$ is a projective morphism, then
\[
\RHom_Y(\F,\fsh G)=\RHom_R(\RG \F,G)
\]
for all $\F\in \D(\Qcoh Y)$, and $G\in\D^+(\Mod R)$.
\end{thm}
\begin{proof}
The sheafified duality theorem \cite[\S6, 6.3]{Neeman2} reads
\[
\Rf\RsHom_Y(\F,\fsh G)=\RsHom_R(\Rf \F,G).
\]
In this situation, $\Rf=\RG$, and since $\Spec R$ is affine, global and local hom agree.  
\end{proof}

The following is also well-known  \cite[3.2.9]{VdB1d}.
\begin{lemma}\label{fshO}
Suppose that $f\colon Y\to\Spec R$ is a projective birational morphism, where $Y$ and $R$ are both Gorenstein varieties of dimension $d$.  Then\\
\t{(1)} $\fsh\O_R$ is a line bundle.\\
\t{(2)} $f$ is crepant if and only if $\fsh\c{O}_R=\c{O}_Y$.
\end{lemma}
\begin{proof}
Since $R$ is Gorenstein $\omega_R$ is a line bundle and thus is a compact object in $\D(\Mod R)$.  Hence by \cite[p227--228]{Neeman2} we have $\fsh\omega_R=\Lf \omega_R\otimes^{\bf L}_Y\fsh\O_R=f^*\omega_R\otimes_Y\fsh\O_R$ and so
\[
\omega_Y=\cd_Y[-d]=\fsh\cd_R[-d]=\fsh\omega_R=f^*\omega_R\otimes_Y\fsh\O_R.
\]
Since both $\omega_Y$ and $f^*\omega_R$ are line bundles, $\fsh\O_R=(f^*\omega_R)^{-1}\otimes_Y\omega_Y$ is a line bundle. Moreover $f$ is crepant if and only if $f^*\omega_R=\omega_Y$ if and only if $\fsh\c{O}_R=\c{O}_Y$.
\end{proof}

The following result show that crepancy implies that $\Lambda$ is Cohen-Macaulay.

\begin{lemma}\label{EndisCMsheaf}
Suppose that $f\colon Y\to\Spec R$ is a crepant projective birational morphism, where $Y$ and $R$ are both Gorenstein varieties of dimension $d$. Then $\End_{\D(Y)}(\c{V})\in\CM R$ for any pretilting complex $\c{V}$ of $Y$.
\end{lemma}
\begin{proof}
We have
\[
\RsHom_Y(\c{V},\c{V})\cong\RsHom_Y(\c{V},\c{V}\otimes^{\bf L}_Y\O_Y)\cong\RsHom_Y(\RsHom_Y(\c{V},\c{V}),\O_Y).
\]
Applying $\RG$ to both sides, we have
\begin{eqnarray*}
\RHom_Y(\c{V},\c{V})&\cong&\RHom_Y(\RsHom_Y(\c{V},\c{V}),\O_Y)\\
&\stackrel{\mbox{\scriptsize\ref{fshO}}}{\cong}&\RHom_Y(\RsHom_Y(\c{V},\c{V}),\fsh\O_R)\\
&\stackrel{\mbox{\scriptsize\ref{SDT}}}{\cong}&\RHom_R(\RHom_Y(\c{V},\c{V}),R).
\end{eqnarray*}
Since $\c{V}$ is pretilting $\RHom_Y(\c{V},\c{V})=\End_{\D(Y)}(\c{V})$, so the above isomorphism reduces to
\[
\End_{\D(Y)}(\c{V})=\RHom_R(\End_{\D(Y)}(\c{V}),R).
\]
Hence applying $H^i$ to both sides, we obtain $\Ext^i_R(\End_{\D(Y)}(\c{V}),R)=0$ for all $i>0$. Thus $\End_{\D(Y)}(\c{V})\in\CM R$, as required.
\end{proof}

To show the converse of \ref{EndisCMsheaf}, namely $\Lambda$ is Cohen--Macaulay implies crepancy, will involve Serre functors.  However, since we are in the singular setting these are somewhat more subtle than usual.  The following is based on \cite[7.2.6]{Ginz}.

\begin{defin}
Suppose that $Y\to\Spec R$ is a morphism where $R$ is CM ring with a canonical module $C_R$.  We say that a functor $\mathbb{S}\colon \Perf(Y)\to\Perf(Y)$ is a \emph{Serre functor relative to $C_R$} if there are functorial isomorphisms
\[
\RHom_R(\RHom_Y(\F,\G),C_R)\cong \RHom_Y(\G,\mathbb{S}(\F))
\]
in $\D(\Mod R)$ for all $\F,\G\in\Perf(Y)$.  If $\Lambda$ is a module-finite $R$-algebra, we define Serre functor $\mathbb{S}\colon\Kb(\proj\Lambda)\to\Kb(\proj\Lambda)$ relative to $C_R$ in a similar way.
\end{defin}

\begin{remark}
We remark that since we are working in the non-local CM setting, canonical modules are not unique.   Although crepancy is defined with respect to the geometric canonical $\omega_R$, and \cite{Ginz} defines Serre functors with respect to $\omega_R$, there are benefits of allowing the flexibility of different canonical modules.  Bridgeland--King--Reid \cite[3.2]{BKR} had technical problems when the geometric canonical is not trivial, since they first had to work locally (where the canonical is trivial), and then extend this globally.  Below, we avoid this problem by considering Serre functors with respect to the canonical module $R$, using the trick in \ref{CMgivescrepancy}(2). 
\end{remark}

The following two lemmas are standard.

\begin{lemma}\label{uniqueSerre}
Suppose that $\mathbb{S}$ and $\mathbb{T}$ are two Serre functors relative to the same canonical $C_R$.  Then $\mathbb{S}$ and $\mathbb{T}$ are isomorphic.
\end{lemma}
\begin{proof}
There are functorial isomorphisms
\[
\RHom_Y(\G,\mathbb{S}(\F))\cong\RHom_Y(\G,\mathbb{T}(\F))
\]
for all $\F,\G\in\Perf(Y)$, which after applying $H^0$ give functorial isomorphisms
\[
\Hom_{\Perf(Y)}(-,\mathbb{S}(-))\cong\Hom_{\Perf(Y)}(-,\mathbb{T}(-))
\]
Since $\mathbb{S}$ and $\mathbb{T}$ take values in $\Perf(Y)$, we may use Yoneda's lemma to conclude that $\mathbb{S}$ and $\mathbb{T}$ are isomorphic.
\end{proof}

\begin{lemma}\label{Serretransfer}
Suppose that $Y\to\Spec R$ is a projective birational morphism between varieties, and $Y$ is  derived equivalent to $\Lambda$. Then any Serre functor $\mathbb{S}\colon\Kb(\proj\Lambda)\to\Kb(\proj\Lambda)$ relative to $C_R$ induces a Serre functor $\mathbb{S}'\colon\Perf(Y)\to\Perf(Y)$ relative to $C_R$
\end{lemma}
\begin{proof}
It is enough to show that we have a triangle equivalence $F\colon\Perf(\Lambda)\to\Perf(Y)$ with a functorial isomorphism $\RHom_\Lambda(A,B)\cong\RHom_Y(FA,FB)$ in $\D(R)$ for all $A,B\in\Perf(\Lambda)$, since then for a
quasi-inverse $E$ of $F$ we have that $\mathbb{S}':=F\circ \mathbb{S}\circ E\colon\Perf(Y)\to \Perf(Y)$ enjoys functorial isomorphisms
\begin{eqnarray*}
\RHom_R(\RHom_Y(FA,FB),C_R)&\cong&
\RHom_R(\RHom_\Lambda(A,B),C_R)\\
&\cong& \RHom_\Lambda(B,\mathbb{S}A)\\
&\cong& \RHom_Y(FB,F(\mathbb{S}A))\\
&\cong& \RHom_Y(FB,\mathbb{S}'(FA))
\end{eqnarray*}
in $\D(R)$.   Let $V$ be a tilting complex of $Y$, and let $\A$ be the DG endomorphism $R$-algebra of $V$.\\
(i) Let $\rC_{\dg}(\A)$ be the DG category of DG $\A$-modules.
We denote by $\pretr(\A)$ the smallest DG subcategory of $\rC_{\dg}(\A)$
which is closed under $[\pm1]$ and cones and contains $\A$.
Since the DG category $\Perf_{\dg}(Y)$ of perfect complexes of $Y$ is
pretriangulated, there exists a fully faithful DG functor
$V\otimes_\A-\colon\pretr(\A)\to\Perf_{\dg}(Y)$
which induces a triangle equivalence $G\colon\Perf(\A)\to\Perf(Y)$
\cite[Section 4.5]{keller}.
In particular, we have a functorial isomorphism
$\Hom^{\bullet}_\A(A,B)\cong\Hom^{\bullet}_Y(V\otimes_\A A,V\otimes_\A B)$ of
DG $R$-modules for all $A,B\in\pretr(\A)$, where we denote by
$\Hom^{\bullet}_\A$ and $\Hom^{\bullet}_Y$ the Hom-sets in our DG categories.
Thus we have a functorial isomorphism
$\RHom_\A(A,B)\cong\RHom_Y(GA,GB)$ in $\D(R)$ for all $A,B\in\Perf(\A)$.\\
(ii) Let $f:\A\to\B$ be a quasi-isomorphism of DG $R$-algebras.
Then the DG functor $\B\otimes_\A-\colon\pretr(\A)\to\pretr(\B)$ gives a triangle
equivalence $H\colon\Perf(\A)\to\Perf(\B)$.
For all $A,B\in\pretr(\A)$, we have a quasi-isomorphism 
$\Hom^{\bullet}_\A(A,B)\to\Hom^{\bullet}_\B(\B\otimes_\A A,\B\otimes_\A B)$ of DG $R$-modules.
In particular, we have a functorial isomorphism
$\RHom_\A(A,B)\cong\RHom_\B(HA,HB)$ in $\D(R)$ for all $A,B\in\Perf(A)$.\\
(iii) Since $V$ is a tilting complex,
we have quasi-isomorphisms $\A^{\le0}\to\A$ and $\A^{\le0}\to\Lambda$
of DG $R$-algebras, where $\A^{\le0}$ is a DG sub $R$-algebra
$(\cdots\to\A^{-1}\to\Ker d^0\to0\to\cdots)$ of $\A$. 
Combining (i) and (ii), we have the desired assertion.
\end{proof}

The following is our key observation.
\begin{lemma}\label{CMgivescrepancy}
Suppose that $f\colon Y\to\Spec R$ is a projective birational morphism, where $Y$ and $R$ are both Gorenstein varieties.\\
\t{(1)} $-\otimes_Y\fsh\O_R\colon \Perf(Y)\to\Perf(Y)$ is a Serre functor relative to $R$.\\
\t{(2)} $f$ is crepant if and only if ${\rm id}\colon \Perf(Y)\to\Perf(Y)$ is a Serre functor relative to $R$.
\end{lemma}
\begin{proof}
(1) We know by \ref{fshO} that $\fsh\O_R$ is a line bundle, so it follows that tensoring by $\fsh\O_R$ gives a functor $-\otimes_Y\fsh\O_R\colon \Perf(Y)\to\Perf(Y)$.  Further, we have functorial isomorphisms
\[
\RHom_Y(\G,\F\otimes_Y\fsh\O_R)\cong\RHom_Y(\RsHom_Y(\F,\G),\fsh\O_R)\stackrel{\mbox{\scriptsize\ref{SDT}}}{\cong}\RHom_R(\RHom_Y(\F,\G),\O_R) 
\]
in $\D(\Mod R)$ for all $\F,\G\in\Perf(Y)$. Thus $-\otimes_Y\fsh\O_R$ is a Serre functor relative to $R$.\\
(2)  By (1) and \ref{uniqueSerre}, ${\rm id}$ is a Serre functor relative to $R$ if and only if $-\otimes_Y\fsh\O_R={\rm id}$ as functors $\Perf(Y)\to\Perf(Y)$. This is equivalent to $\fsh\O_R=\O_Y$, which is equivalent to that $f$ is crepant by \ref{fshO}.
\end{proof}

The following explains the geometric origin of the definition of modifying modules, and is the first main result of this subsection.

\begin{thm}\label{CMiffcrepant}
Let $f\colon Y\to\Spec R$ be a projective birational morphism between $d$-dimensional Gorenstein varieties.  Suppose that $Y$ is derived equivalent to some ring $\Lambda$, then the following are equivalent.\\
\t{(1)} $f$ is crepant.\\
\t{(2)} $\Lambda\in\CM R$.\\
\t{(3)} ${\rm id}\colon \Perf(Y)\to\Perf(Y)$ is a Serre functor relative to $R$.\\
In this case $\Lambda\cong\End_R(M)$ for some $M\in\refl R$.
\end{thm}
\begin{proof}
(3)$\Leftrightarrow$(1) is shown in \ref{CMgivescrepancy}(2), and (1)$\Rightarrow$(2) is shown in \ref{EndisCMsheaf}.\\
(2)$\Rightarrow$(3) Let
\begin{eqnarray*}
\mathbb{S}:=\RHom_R(\Lambda,R)\otimes^{\mathbf{L}}_\Lambda-\colon \D^-(\mod\Lambda)\to\D^-(\mod\Lambda).
\end{eqnarray*}
By \cite[3.5(2)(3)]{IR}, there exists a functorial isomorphism
\begin{eqnarray}
\RHom_\Lambda(A,\mathbb{S}(B))\cong \RHom_R(\RHom_\Lambda(B,A),R)\label{IR3.5}
\end{eqnarray}
in $\D(R)$ for all $A\in\Db(\mod\Lambda)$ and all $B\in\Kb(\proj\Lambda)$.  

Now suppose that $\Lambda\in\CM R$. Since it is reflexive, $\Lambda\cong\End_R(M)$ for some $M\in\refl R$ by \ref{Endincorrectform}. Then we have
\[
\Lambda\cong\Hom_R(\Lambda,R)\cong\RHom_R(\Lambda,R)\label{symm iso}
\]
in $\D(\Lambda\otimes_R\Lambda^{\op})$ where the first isomorphism holds by \ref{ModisdCY}(1), and the second since $\Lambda\in\CM R$. Thus we have $\mathbb{S}\cong{\rm id}$, and \eqref{IR3.5} shows that ${\rm id}\colon\Kb(\proj\Lambda)\to\Kb(\proj\Lambda)$ is a Serre functor relative to $R$. By \ref{Serretransfer} it follows that ${\rm id}\colon\Perf(Y)\to\Perf(Y)$ is a Serre functor relative to $R$.
\end{proof}

When further $Y$ is smooth, we can strengthen \ref{CMiffcrepant}.

\begin{cor}\label{NCCRdimd}
Let $f\colon Y\to\Spec R$ be a projective birational morphism between $d$-dimensional Gorenstein varieties.  Suppose that $Y$ is derived equivalent to some ring $\Lambda$, then the following are equivalent.\\
\t{(1)} $f$ is a crepant resolution of $\Spec R$.\\
\t{(2)} $\Lambda$ is an NCCR of $R$.
\end{cor}
\begin{proof}
By \ref{CMiffcrepant} $f$ is crepant if and only if $\Lambda\in \CM R$. Assume that this is satisfied, then again by \ref{CMiffcrepant} $\Lambda\cong\End_R(M)$ for some $M\in\refl R$.  By \ref{Dbisolated}(1) $Y$ is smooth if and only if $\Lambda$ is non-singular, which means that $\Lambda$ is an NCCR of $R$.
\end{proof}

We are now in a position to relate maximal modification algebras and $\mathds{Q}$-factorial terminalizations in dimension three.  The following is the second main result of this subsection.

\begin{thm}\label{mainQfact}
Let $f\colon Y\to\Spec R$ be a projective birational morphism, where $Y$ and $R$ are both Gorenstein varieties of dimension three. Assume that $Y$ has (at worst) isolated singularities $\{ x_1,\hdots,x_n\}$ where each $\c{O}_{Y,x_i}$ is a hypersurface.  If $Y$ is derived equivalent to some ring $\Lambda$, then the following are equivalent.\\
\t{(1)} $f$ is crepant and $Y$ is $\mathds{Q}$-factorial.\\
\t{(2)} $\Lambda$ is an MMA of $R$.\\
In this situation, all MMAs of $R$ are derived equivalent and have isolated singularities.
\end{thm}
\begin{proof} 
By \ref{CMiffcrepant} $f$ is crepant if and only if $\Lambda\in \CM R$. Assume that this is satisfied, then again by \ref{CMiffcrepant} $\Lambda\cong\End_R(M)$ for some $M\in\refl R$.  Thus $\Dsg(Y)\simeq\uCM\Lambda$ by \ref{DbinducesPerf}, \ref{ModisdCY}(2) and \ref{Buchweitz}.  Since the singularities of $Y$ are isolated, $\Lambda$ has isolated singularities by \ref{Dbisolated}(2). Now by \ref{Section3_main}(1), $Y$ is $\mathds{Q}$-factorial if and only if $\Dsg(Y)\simeq\uCM\Lambda$ is rigid-free, which by \ref{MMtestforuCM} holds if and only if $\Lambda$ is an MMA.   Lastly, all maximal modification algebras are derived equivalent \cite[4.15]{IW}, thus all are derived equivalent to $Y$ and hence all have isolated singularities by \ref{Dbisolated}(2).
\end{proof}

\begin{remark}\label{fieldC}
When $\k=\mathbb{C}$ and $Y$ has only terminal singularities, the geometric assumptions in \ref{mainQfact} are satisfied since terminal singularities are isolated for 3-folds \cite[5.18]{KM}, and Gorenstein terminal singularities are Zariski locally hypersurfaces \cite[0.6(I)]{Pagoda}.  This shows that when $\k=\mathbb{C}$ and $Y$ has only terminal singularities, $Y$ is a $\mathds{Q}$-factorial terminalization of $\Spec R$ if and only if $\Lambda$ is an MMA.
\end{remark}

\begin{remark}
We also remark that if Conjecture \ref{3foldconj} is true, then the $Y$ in \ref{mainQfact} always admit a tilting complex.   Note that the corresponding statement of \ref{3foldconj} in dimension four is false, see \ref{counterexample}.
\end{remark} 

Note that all the results in this subsection remain valid
in the complete local setting. In particular \ref{corforSection4} below corresponds to \ref{mainQfact}, and will be used in \S\ref{Dbsection}.  Here, a morphism $f\colon Y\to X=\Spec R$
is called \emph{crepant} if $f^*\omega_X=\omega_Y$ holds (as before), where now $X$ has a unique canonical sheaf $\omega_X=\mathcal{O}_R$, and we choose the canonical sheaf $\omega_Y:=f^!\omega_X$ on $Y$.
 
\begin{cor}\label{corforSection4}
Suppose $f\colon Y\to\Spec R$ is a projective birational morphism between $3$-dimensional Gorenstein normal integral schemes, satisfying (ELF), where $R$ is a complete local ring containing a copy of its residue field. Suppose further that $Y$ has (at worst) isolated singularities $\{ x_1,\hdots,x_n\}$ where each $\c{O}_{Y,x_i}$ is a hypersurface.  If $Y$ is derived equivalent to some ring $\Lambda$, then the following are equivalent.\\
\t{(1)}  $Y$ is complete locally $\mathds{Q}$-factorial and $f$ is crepant.\\
\t{(2)}  $Y$ is $\mathds{Q}$-factorial and $f$ is crepant.\\
\t{(3)}  $\Lambda$ is an MMA of $R$.\\
In this situation, all MMAs of $R$ are derived equivalent and have isolated singularities.
\end{cor}
\begin{proof}
(2)$\Leftrightarrow$(3) is similar to  \ref{mainQfact} (1)$\Leftrightarrow$(2).\\
(1)$\Leftrightarrow$(2) By \ref{Section3_main}(1)(2), we only have to show that $\Dsg(Y)=\overline{\Dsg(Y)}$ holds.  By \ref{CMiffcrepant} $\Lambda\cong\End_R(M)\in\CM R$ for some $M\in\refl R$. We have $\Dsg(Y)\simeq \uCM\Lambda$ by \ref{DbinducesPerf}, \ref{ModisdCY}(2) and \ref{Buchweitz}.
The endomorphism ring $\End_\Lambda(X)$ of any $X\in\mod\Lambda$ is again a module-finite algebra over a complete local ring $R$. Thus any idempotent in $\underline{\End}_\Lambda(X)$ is an image of some idempotent in $\End_\Lambda(X)$ by \cite[6.5, 6.7]{CR}, which corresponds to some direct summand of $X$. Therefore the category $\uCM\Lambda\simeq\Dsg(Y)$ is idempotent complete, and we have the assertion.
\end{proof}

\subsection{A counterexample in dimension four}

Here we show that we cannot always expect the setup of \ref{mainQfact} and \ref{NCCRdimd} to hold in higher dimension.  This puts severe limitations on any general homological theory that covers dimension four.  

The following is an extension of an example of Dao \cite[3.5]{DaoNCCR}.

\begin{thm}\label{counterexample}
Let $R:=\C{}[x_0,x_1,x_2,x_3,x_4]/(x_0^5+x_1^4+x_2^4+x_3^4+x_4^4)$.  This has a crepant resolution, which we denote by $Y\to\Spec R$.  Then there is no algebra $\Lambda$ that is derived equivalent to $Y$.
\end{thm}
\begin{proof}
Note first that $R$ is an isolated singularity, with unique singular point at the origin.  A projective crepant resolution exists by \cite[Thm.\ A.4]{Lin}.  Suppose $\Lambda$ is derived equivalent to $Y$, then by \ref{NCCRdimd} $\Lambda\cong\End_R(M)$ is an NCCR of $R$.  Completing with respect to the maximal ideal of the origin, this implies that $\End_{\widehat{R}}(\widehat{M})$ is an  NCCR of $\widehat{R}$.  But $\End_{\widehat{R}}(\widehat{M})\in\CM \widehat{R}$ with $\widehat{M}\in\refl \widehat{R}$, so since $\widehat{R}$ is a 4-dimensional local isolated hypersurface, $\widehat{M}$ must be free \cite[2.7(3)]{DaoNCCR}. But this forces $\End_{\widehat{R}}(\widehat{M})\cong M_n(\widehat{R})$, which is a contradiction since $M_n(\widehat{R})$ has infinite global dimension.
\end{proof}

\section{$cA_n$ singularities via derived categories}\label{Dbsection}

We now illustrate the results of the previous section to give many examples of maximal modification algebras. The following result is due to Shepherd--Barron (unpublished).

\begin{prop}\label{SB}
Let $R=\C{}[[u,v,x,y]]/(uv-f(x,y))$ be an isolated $cA_{n}$ singularity.  Then the following are equivalent\\
\t{(1)}  There does not exist a non-trivial crepant morphism $Y\to\Spec R$.\\
\t{(2)} $f(x,y)$ is irreducible.\\
\t{(3)} $R$ is factorial.\\
\t{(4)} $R$ is $\mathds{Q}$-factorial.\\
\t{(5)} $R$ is an MM $R$-module.
\end{prop}
\begin{proof}
(1)$\Rightarrow$(2) If $f(x,y)$ factors non-trivially as $f=f_{1}f_{2}$ then blowing up the ideal $(u,f_{1})$ yields a non-trivial crepant morphism (as in the calculation in \S\ref{crepmod}).  Hence $f(x,y)$ cannot factor.\\
(2)$\Rightarrow$(3) This is easy (see e.g.\ proof of \cite[5.9]{IW6}). \\
(3)$\Rightarrow$(4)  This is clear.\\
(4)$\Rightarrow$(5) This follows from \ref{UFDandrigid}(5)$\Rightarrow$(1).\\
(5)$\Rightarrow$(1) Since $R$ is an isolated cDV singularity of type $A$ (see e.g. \cite[6.1(e)]{BIKR}), it is a terminal singularity by \cite[1.1]{Pagoda}. If there exists a non-trivial crepant morphism $f\colon Y\to\Spec R$, then $f$ must have one-dimensional fibres.  Hence by \cite[3.2.10]{VdB1d} there exists a non-projective modifying module.  Since $R$ is a MM $R$-module, this contradicts \ref{RisMMthenMMareFree}.
\end{proof}

\subsection{Crepant modifications of $cA_n$ singularities}\label{crepmod}

In this section we work over $\KK$, an algebraically closed field of characteristic zero. Let
$$R:=\KK[[u,v,x,y]]/(uv-f(x,y)),$$
with $f\in\m:=(x,y)\subseteq \KK[[x,y]]$.  We let $f=f_{1}\hdots f_{n}$ be a factorization into prime elements of $\KK[[x,y]]$.  We restrict to this complete local setting since it simplifies the proof of the main theorem \ref{Dbforflags}.  A similar version of \ref{Dbforflags} is true when $R$ is not complete local, though this requires a much more complicated proof \cite{W}.

We remark that when $\KK=\mathbb{C}$, $R$ is a $cA_m$ singularity for $m:={\rm ord}(f)-1$ (see e.g. \cite[6.1(e)]{BIKR}). This is terminal if and only if $R$ is an isolated singularity \cite[1.1]{Pagoda} if and only if $(f_i)\neq(f_j)$ for all $i\neq j$.

For any subset $I\subseteq\{ 1,\hdots, n\}$ we denote
\[
f_I:=\prod_{i\in I}f_i\ \mbox{ and }\ T_I:=(u,f_I)
\]
which is an ideal of $R$.  For a collection of subsets 
\[
\emptyset\subsetneq I_1\subsetneq I_2\subsetneq...\subsetneq
I_m\subsetneq\{1,2,...,n\},
\]
we say that $\F=(I_1,\hdots,I_m)$ is a {\em flag in the set $\{ 1,2,\hdots, n\}$}. We say that the flag $\F$ is {\em maximal} if $n=m+1$.  Given a flag $\F=(I_1,\hdots,I_m)$, we define
\[
T^\c{F}:=R\oplus\left(\bigoplus_{j=1}^{m} T_{I_j}\right) .
\]
On the other hand, geometrically, given a flag $\F=(I_1,\hdots,I_m)$ we define a scheme $X^\F$ as follows:

First, let $X^{\F_1}\to\Spec R$ be the blowup of the ideal $(u,f_{I_1})$ in $R=\KK[[u,v,x,y]]/(uv-f)$.  Then $X^{\F_1}$ is covered by two open charts, given explicitly by 
\[
R_1:=\frac{\KK[[u,v,x,y]][V_1]}{\left(\begin{array}{l}v-V_1\frac{f}{f_{I_1}}\\uV_1-f_{I_1}\end{array}\right)} 
\quad\text{and}\quad 
R_2:=\frac{\KK[[u,v,x,y]][U_1]}{\left(\begin{array}{l}u-f_{I_1}U_1\\ U_1v-\frac{f}{f_{I_1}}\end{array}\right)}.
\]
The new coordinates $V_1=\frac{f_{I_1}}{u}$ and $U_1=\frac{u}{f_{I_1}}$ glue to give a copy of $\mathbb{P}^1$ inside $X^{\F_1}$ (possibly moving in a family), which maps to the origin of $\Spec R$.  

Next, let $X^{\F_2}\to X^{\F_1}$ be the blowup of the divisor $(U_1,\frac{f_{I_2}}{f_{I_1}})$ in the second co-ordinate chart $R_2$.  Note that the zero set of the divisor $(U_1,\frac{f_{I_2}}{f_{I_1}})$ does not intersect the first co-ordinate chart $R_1$, so $R_1$ is unaffected under the second blowup.  Locally above $R_2$, the calculation to determine the structure of $X^{\F_2}$ is similar to the above. Thus $X^{\F_2}$ is covered by the three affine open sets 
\[
\frac{\KK[[u,v,x,y]][V_1]}{\left(\begin{array}{l}v-V_1\frac{f}{f_{I_1}}\\uV_1-f_{I_1}\end{array}\right)} 
\quad\text{and}\quad
\frac{\KK[[u,v,x,y]][U_1,V_2]}{\left(\begin{array}{l}u-f_{I_1}U_1\\ v-\frac{f}{f_{I_2}}V_2\\ U_1V_2-\frac{f_{I_2}}{f_{I_1}}\end{array}\right)} 
\quad\text{and}\quad
\frac{\KK[[u,v,x,y]][U_2]}{\left(\begin{array}{l}u-f_{I_2}U_2\\ U_2v-\frac{f}{f_{I_2}}\end{array}\right)}
\]
The $U_2$ and the $V_2$ coordinates again glue to produce another $\mathbb{P}^1$ (again which might move in a family) which maps to the origin of $R_2$ and hence to the origin of $R$.

Continuing by blowing up the ideal $(U_2,\frac{f_{I_3}}{f_{I_2}})$, in this way we obtain a chain of projective birational morphisms
\[
X^{\F_m}\to X^{\F_{m-1}}\to\hdots\to X^{\F_1}\to \Spec R,
\]
and we define $X^\F:=X^{\F_m}$. See \ref{pictures} later for a picture of this process. 

Note that $X^{\F}$, being projective over the base $\Spec R$, automatically satisfies (ELF) \cite[2.1.3]{TT}, and by inspection of the charts, $X^\F$ is a normal integral Gorenstein scheme, so we can apply the results in the previous section.  Also note that if we complete each of the above affine open sets at the origin we obtain 
\[
\frac{\KK[[u,V_1,x,y]]}{\left(uV_1-f_{I_1}\right)}, \quad
\frac{\KK[[U_1,V_2,x,y]]}{(U_1V_2-\frac{f_{I_2}}{f_{I_1}})},\quad\hdots\quad
\frac{\KK[[U_{m-1},V_m,x,y]]}{(U_{m-1}V_m-\frac{f_{I_m}}{f_{I_{m-1}}})},\quad
\frac{\KK[[U_m,v,x,y]]}{(U_mv-\frac{f}{f_{I_m}})}.
\]

\begin{thm}\label{Dbforflags}
Given a flag $\F=(I_1,\hdots,I_m)$, denote $X^\F$ and $T^\F$ as above.  Then $X^\F$ is derived equivalent to $\End_R(T^\F)$.
\end{thm}
\begin{proof}
In the explicit calculation for $X^{\F}$ above, the preimage of the unique closed point $\n$ of $\Spec R$ is a chain of $\mathbb{P}^1$'s (some of which move in a family), in a type $A_m$ configuration.  Now by \cite[Thm.\ B]{VdB1d}, there is a tilting bundle on $X^{\F}$ given as follows:  let $C=f^{-1}(\n)$.  Giving $C$ the reduced scheme structure, write $C_{\rm red}=\bigcup_{i\in I}C_i$, and let $\c{L}_i$ denote the line bundle on $X^{\F}$ such that $\c{L}_i\cdot C_j=\delta_{ij}$.  If the multiplicity of $C_i$ in $C$ is equal to one, set $\c{M}_i:=\c{L}_i$ \cite[3.5.4]{VdB1d}, else define $\c{M}_i$ to be given by the extension
\[
0\to\c{O}^{r_i-1}\to\c{M}_i\to\c{L}_i\to 0
\]
associated to a minimal set of $r_i-1$ generators of $H^1(X^{\F},\c{L}_i^{-1})$.  Then $\c{O}\oplus(\bigoplus_{i\in I}\c{M}_i)$ is a tilting bundle on $X^\F$ \cite[3.5.5]{VdB1d}.

In our situation, let $C_i$ be the curve in $X^\F$, above the origin of $\Spec R$, which in the process of blowing up first appears in $X^{\F_i}$ (see \ref{pictures} for a picture of this). We now claim that all the curves $C_i$ have multiplicity one.  If some $C_i$ had multiplicity greater than one, then $\c{M}_i$ would be an indecomposable bundle with rank greater than one \cite[3.5.4]{VdB1d}.  But this would imply, by \cite[3.2.9]{VdB1d}, that its global sections $H^0(\c{M}_i):=M_i$ is an indecomposable CM $R$-module of rank greater than one, such that $\End_R(M_i)\in \CM R$.  But this is impossible (see e.g.\ \cite[5.24]{IW6}).  Hence all curves have multiplicity one, and so $\c{O}\oplus(\bigoplus_{i=1}^m \c{L}_i)$ is a tilting bundle on $X^\F$.

Now $\End_{X^\F}(\c{O}\oplus(\bigoplus_{i=1}^m \c{L}_i))\cong \End_R(R\oplus(\bigoplus_{i=1}^m H^0(\c{L}_i)) )$, and so it remains to show that $R\oplus(\bigoplus_{i=1}^m H^0(\c{L}_i))\cong T^\F$.  In fact, we claim that $H^0(\c{L}_i)\cong (u,f_{I_i})$.  But it is easy to see (for example using the \v{C}ech complex) that the rank one CM module $H^0(\c{L}_i)$ is generated as an $R$-module by $u$ and $f_{I_i}$. 
\end{proof}

So as to match our notation with \cite{BIKR} and \cite{DH}, we can (and do) identify maximal flags with elements of the symmetric group $\mathfrak{S}_n$. Hence we regard each $\omega\in\mathfrak{S}_n$ as the maximal flag 
\[
\{\omega(1) \}\subset\{\omega(1),\omega(2) \}\subset\hdots\subset\{\omega(1),\hdots,\omega(n-1) \}.
\]

The following is simply a special case of \ref{Dbforflags}.
\begin{cor}\label{keyDb}
The scheme $X^{\omega}$ is derived equivalent to $\End_R(T^\omega)$ for all $\omega\in\mathfrak{S}_n$. Moreover the completions of singular points of $X^{\omega}$ are precisely $\KK[[u,v,x,y]]/(uv-f_i)$ for all $1\leq i\leq n$ such that $f_i\in(x,y)^2$.
\end{cor}

This yields the following generalization of \cite[1.5]{BIKR} and \cite[4.2]{DH}, and provides a conceptual reason for the condition $f_i\notin\m^2$.  
\begin{cor}\label{derivedMMCTforcAn}
Let $f_1,\hdots,f_n\in\m:=(x,y)\subseteq \KK[[x,y]]$ be irreducible polynomials and $R=\KK[[u,v,x,y]]/(uv-f_1\hdots f_n)$. Then\\
\t{(1)}  Each $T^\omega$ ($\omega\in\mathfrak{S}_n$) is an MM $R$-module which is a generator.  The endomorphism rings $\End_R(T^\omega)$ have isolated singularities.\\
\t{(2)}  $T^\omega$ is a CT $R$-module for some $\omega\in\mathfrak{S}_n$ $\iff$ $T^\omega$ is a CT $R$-module for all $\omega\in\mathfrak{S}_n$ $\iff$ $f_i\notin\m^2$ for all $1\leq i\leq n$.
\end{cor}
\begin{proof}
(1) $\End_R(T^\omega)$ is derived equivalent to a scheme whose (complete local) singularities are listed in \ref{keyDb}.  Since each $f_i$ is irreducible, these are all factorial by \ref{SB}.  Hence the result follows from \ref{corforSection4}.\\
(2) $T^\omega$ is CT if and only if $\Gamma:=\End_R(T^\omega)$ is an NCCR \cite[5.4]{IW}. By \ref{Dbisolated}(1) this occurs if and only if $X^\omega$ is smooth.  But by \ref{keyDb}, this happens if and only if each $uv=f_i$ is smooth, which is if and only if each $f_i\notin \m^2$.
\end{proof}

For an algebraic non-derived category proof of the above, see \cite[\S5]{IW6}.  In fact  the commutative algebraic method in \cite{IW6} is stronger, since it also gives that there are no other MM (respectively CT) generators. Also, the methods in \cite{IW6} work over more general fields.

\begin{remark}\label{pictures} It is useful to visualize the above in an example.   

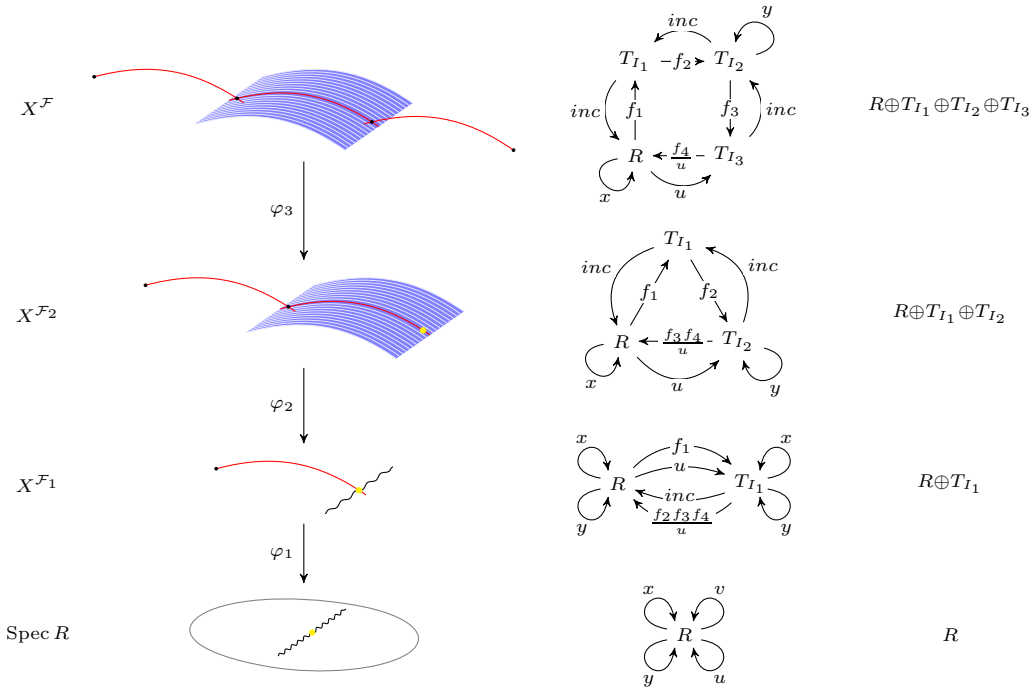
\begin{center}
\begin{tikzpicture}
\node (a2) at (0,7) 
{\begin{tikzpicture} [transform shape, rotate=-10]
\fill[fill=blue!50] (0,0,1) -- (0,0,-1) to [bend left=25] (2,0,-1) -- (2,0,1) to [bend right=25] (0,0,1);
\foreach \y in {0.1,0.2,...,1}{ 
\draw[very thin,blue!10] (0,0,\y) to [bend left=25] (2,0,\y);
\draw[very thin,blue!10] (0,0,-\y) to [bend left=25] (2,0,-\y);}
\draw[red] (0,0,0) to [bend left=25] (2,0,0);
\draw[red] (1.8,0,0) to [bend left=25] (3.8,0,0);
\draw[red] (-1.8,0,0) to [bend left=25] (0.2,0,0);
\filldraw [black] (0.1,0.045,0) circle (0.5pt);
\filldraw [black] (1.9,0.045,0) circle (0.5pt);
\filldraw [black] (-1.8,0,0) circle (0.5pt);
\filldraw [black] (3.8,0,0) circle (0.5pt);
\end{tikzpicture} };

\node (b2) at (0,4.25)
{\begin{tikzpicture} [transform shape, rotate=-10]
\fill[fill=blue!50] (0,0,1) -- (0,0,-1) to [bend left=25] (2,0,-1) -- (2,0,1) to [bend right=25] (0,0,1);
\foreach \y in {0.1,0.2,...,1}{ 
\draw[very thin,blue!10] (0,0,\y) to [bend left=25] (2,0,\y);
\draw[very thin,blue!10] (0,0,-\y) to [bend left=25] (2,0,-\y);}
\draw[red] (0,0,0) to [bend left=25] (2,0,0);
\draw[red] (-1.8,0,0) to [bend left=25] (0.2,0,0);
\filldraw [black] (0.1,0.045,0) circle (0.5pt);
\filldraw [yellow] (1.9,0.045,0) circle (1pt);
\filldraw [black] (-1.8,0,0) circle (0.5pt);
\end{tikzpicture} };

\node (c2) at (0,2)
{\begin{tikzpicture} [transform shape, rotate=-10]
\draw[red] (-1.8,0,0) to [bend left=25] (0.2,0,0);
\draw[decorate, 
decoration={snake,amplitude=.2mm,segment length=2.1mm}] (0.1,0.045,1) -- node {} (0.1,0.045,-1);
\filldraw [yellow] (0.1,0.045,0) circle (1pt);
\filldraw [black] (-1.8,0,0) circle (0.5pt);
\end{tikzpicture} };

\node (d2) at (0,0)
{\begin{tikzpicture} [transform shape, rotate=-10]
\draw[decorate, 
decoration={snake,amplitude=.2mm,segment length=1mm}] (0.1,0.045,1) -- node {} (0.1,0.045,-1);
\draw[gray] (0.1,0,1.3) .. controls (0.8215,0,1.3) and (1.4,0,0.7215) .. (1.4,0,0);
\draw[gray] (1.4,0,0) .. controls (1.4,0,-0.7215) and (0.8215,0,-1.3) .. (0.1,0,-1.3);
\draw[gray] (0.1,0,-1.3) .. controls (-0.8215,0,-1.3) and (-1.4,0,-0.7215) .. (-1.4,0,0);
\draw[gray] (-1.4,0,0) .. controls (-1.4,0,0.7215) and (-0.8215,0,1.3) .. (0.1,0,1.3);
\filldraw [yellow] (0.1,0.045,0) circle (1pt);
\end{tikzpicture}}; 

\node (a1) at (-3.5,7) {$\scriptstyle X^{\F}$}; 
\node (b1) at (-3.5,4.25) {$\scriptstyle X^{\F_2}$}; 
\node (c1) at (-3.5,2) {$\scriptstyle X^{\F_1}$}; 
\node (d1) at (-3.5,0) {$\scriptstyle \Spec R$}; 

\draw[->] (a2) -- node[left] {$\scriptstyle \varphi_{3}$} (b2);
\draw[->] (b2) -- node[left] {$\scriptstyle \varphi_{2}$} (c2);
\draw[->] (c2) -- node[left] {$\scriptstyle \varphi_{1}$} (d2);

\node (a4) at (5,7)  
{\begin{tikzpicture} 
\node[name=s,regular polygon, regular polygon sides=4, minimum size=1.75cm] at (0,0) {}; 
\node (C1) at (s.corner 1)  {$\scriptstyle T_{I_2}$};
\node (C2) at (s.corner 2)  {$\scriptstyle T_{I_1}$};
\node (C4) at (s.corner 4)  {$\scriptstyle T_{I_3}$};
\node (C3) at (s.corner 3)  {$\scriptstyle R$};
\draw[->] (C4) -- node[inner sep=0.5pt,fill=white,pos=0.45] {$\scriptstyle \frac{f_{4}}{u}$} (C3); 
\draw[->] (C3) -- node[inner sep=0.5pt,fill=white,pos=0.5] {$\scriptstyle f_{1}$} (C2); 
\draw[->] (C2) -- node[inner sep=0.5pt,fill=white,pos=0.5] {$\scriptstyle f_{2}$} (C1); 
\draw[->] (C1) -- node[inner sep=0.5pt,fill=white,pos=0.5] {$\scriptstyle f_{3}$} (C4);
\draw [->,bend right=45,looseness=1,above] (C1) to node[inner sep=1.5pt,fill=white,pos=0.5]  {$\scriptstyle inc$} (C2);
\draw [->,bend right=45,looseness=1,left] (C2) to node[inner sep=1.5pt,fill=white,pos=0.5]  {$\scriptstyle inc$} (C3);
\draw [->,bend right=45,looseness=1,below] (C3) to node[inner sep=1.5pt,fill=white,pos=0.5]  {$\scriptstyle u$} (C4);
\draw [->,bend right=45,looseness=1,right] (C4) to node[inner sep=1.5pt,fill=white,pos=0.5]  {$\scriptstyle inc$} (C1);
\draw[->]  (C3) edge [in=-100,out=-170,loop,looseness=6] node[below,pos=0.5] {$\scriptstyle x$} (C3);
\draw[->]  (C1) edge [in=75,out=15,loop,looseness=5] node[above,pos=0.5] {$\scriptstyle y$} (C1);
\end{tikzpicture}};

\node (b4) at (5,4.25)  
{\begin{tikzpicture} 
\node[name=s,regular polygon, regular polygon sides=3, minimum size=1.75cm] at (0,0) {}; 
\node (c1) at (s.corner 1)  {$\scriptstyle T_{I_1}$};
\node (c2) at (s.corner 2)  {$\scriptstyle R$};
\node (c3) at (s.corner 3)  {$\scriptstyle T_{I_2}$};
\draw[->] (c2) -- node  [inner sep=0.5pt,fill=white,pos=0.5] {$\scriptstyle f_{1}$} (c1); 
\draw[->] (c1) -- node  [inner sep=0.5pt,fill=white,pos=0.5] {$\scriptstyle f_{2}$} (c3);
\draw[->] (c3) -- node  [inner sep=0.5pt,fill=white,pos=0.4] {$\scriptstyle \frac{f_{3}f_{4}}{u}$} (c2); 
\draw [->,bend right=45,looseness=1.2] (c1) to node [inner sep=1.5pt,fill=white,pos=0.5,above left] {$\scriptstyle inc$} (c2);
\draw [->,bend right=45,looseness=1.2] (c2) to node [inner sep=1.5pt,fill=white,pos=0.5,below] {$\scriptstyle u$} (c3);
\draw [->,bend right=45,looseness=1.2] (c3) to node [inner sep=1.5pt,fill=white,pos=0.5,above right] {$\scriptstyle { inc}$} (c1);
\draw[->]  (c2) edge [in=-100,out=-170,loop,looseness=6] node[below,pos=0.5] {$\scriptstyle x$} (c2);
\draw[->]  (c3) edge [in=-75,out=-15,loop,looseness=5] node[below,pos=0.5] {$\scriptstyle y$} (c3);
\end{tikzpicture}};

\node (c4) at (5,2)  
{\begin{tikzpicture} 
\node (C1) at (0,0)  {$\scriptstyle R$};
\node (C1a) at (-0.1,0.05)  {};
\node (C1b) at (-0.1,-0.05)  {};
\node (C2) at (1.75,0)  {$\scriptstyle T_{I_1}$};
\node (C2a) at (1.85,0.05) {};
\node (C2b) at (1.85,-0.05) {};
\draw [->,bend left=45,looseness=1,pos=0.5] (C1) to node[inner sep=0.5pt,fill=white]  {$\scriptstyle f_{1}$} (C2);
\draw [->,bend left=20,looseness=1,pos=0.5] (C1) to node[inner sep=0.5pt,fill=white]  {$\scriptstyle u$} (C2);
\draw [->,bend left=45,looseness=1,pos=0.5] (C2) to node[inner sep=0.5pt,fill=white]  {$\scriptstyle \frac{f_{2}f_{3}f_{4}}{u}$} (C1);
\draw [->,bend left=20,looseness=1,pos=0.5] (C2) to node[inner sep=0.5pt,fill=white,below=-5pt]  {$\scriptstyle inc$} (C1);
\draw[<-]  (C1b) edge [in=-170,out=-100,loop,looseness=12,pos=0.5] node[below] {$\scriptstyle y$} (C1b);
\draw[<-]  (C1a) edge [in=170,out=100,loop,looseness=12,pos=0.5] node[above] {$\scriptstyle x$} (C1a);
\draw[->]  (C2b) edge [in=-80,out=-10,loop,looseness=12,pos=0.5] node[below] {$\scriptstyle y$} (C2b);
\draw[->]  (C2a) edge [in=80,out=10,loop,looseness=12,pos=0.5] node[above] {$\scriptstyle x$} (C2a);
\end{tikzpicture}};

\node (d4) at (5,0)  
{\begin{tikzpicture} 
\node (1) at (0,0)  {$\scriptstyle R$};
\node (1a) at (-0.1,0.05)  {};
\node (1b) at (-0.1,-0.05)  {};
\node (2a) at (0.1,0.05)  {};
\node (2b) at (0.1,-0.05)  {};
\draw[<-]  (1b) edge [in=-170,out=-100,loop,looseness=12,pos=0.5] node[below] {$\scriptstyle y$} (1b);
\draw[<-]  (1a) edge [in=170,out=100,loop,looseness=12,pos=0.5] node[above] {$\scriptstyle x$} (1a);
\draw[->]  (2b) edge [in=-80,out=-10,loop,looseness=12,pos=0.5] node[below] {$\scriptstyle u$} (2b);
\draw[->]  (2a) edge [in=80,out=10,loop,looseness=12,pos=0.5] node[above] {$\scriptstyle v$} (2a);
\end{tikzpicture}};

\node (a5) at (8.5,7) {$\scriptstyle R\oplus T_{I_1}\oplus T_{I_2}\oplus T_{I_3}$}; 
\node (b5) at (8.5,4.25) {$\scriptstyle R\oplus T_{I_1}\oplus T_{I_2}$}; 
\node (c5) at (8.5,2) {$\scriptstyle R\oplus T_{I_1}$}; 
\node (d5) at (8.5,0) {$\scriptstyle R$}; 

\end{tikzpicture}
\end{center}

Consider $uv=f_{1}f_{2}f_{3}f_{4}$ and choose maximal flag $\F=(I_1:=\{1\}\subset I_2:=\{1,2\}\subset I_3:=\{1,2,3\})$. Thus $T_{I_1}=(u,f_1)$, $T_{I_2}=(u,f_1f_2)$ and $T_{I_3}=(u,f_1f_2f_3)$. Moreover we assume, so that we can draw an accurate picture, that $(f_2)=(f_3)$, and $(f_3)\neq(f_4)$, so that $uv=f_2f_3f_4$ is non-isolated, but $uv=f_3f_4$ is isolated.

On the geometric side we have drawn the non-local picture; the complete local picture is obtained by drawing a tubular neighbourhood around all the red curves, and restricting to the origin in $\Spec R$.

Now on $X^{\F_1}$, the black dot is the singularity $uv=f_1$, whereas the yellow dot is $uv=f_2f_3f_4$.  On $X^{\F_2}$, the latter singularity splits into $uv=f_2$ (middle black dot) and $uv=f_3f_4$ (yellow).  The left hand black dot is still $uv=f_1$. On $X^\F$, the black dots are (reading left to right) $uv=f_1$, $uv=f_2$, $uv=f_3$ and $uv=f_4$.  All yellow dots correspond to non-$\mathds{Q}$-factorial singularities, and all black dots correspond to isolated $\mathds{Q}$-factorial singularities, possibly smooth.  The yellow dots with squiggles through them are non-isolated singularities.

The red curves are the $\mathbb{P}^{1}$'s which map to the origin in $\Spec R$.  The right hand curve in $X^{\F_2}$ moves in a family (represented as blue lines) and $\varphi_{2}$ contracts the whole family; consequently $\varphi_{2}$ contracts a divisor.  However both $\varphi_{1}$ and $\varphi_{3}$ contract a single curve but no divisor and so are flopping contractions.  Note that $\Spec R$ and $X^{\F_1}$ have canonical singularities, but $X^{\F_2}$ and $X^{\F}$ have only terminal singularities.  

The precise form of the quiver relies on the (easy) calculation in \cite[5.33]{IW6}, but this can be ignored for now.  Note that the geometric picture will change depending both on the choice of polynomials, and their ordering.  Blowing up in a different order can change which curves move in families, and also normal bundles of the curves.  Of course, this depends on the choice of the polynomials too.  On the algebraic side, if we change the polynomials or their order then the number of loops change, as do the relations.  Changing the ordering corresponds to mutation (see \cite[5.31]{IW6}).

On each level $i$ the geometric space $X^{\F_i}$ is derived equivalent to the corresponding algebra by \ref{Dbforflags} .  The top space has only $\mathds{Q}$-factorial terminal singularities, hence the top algebra is an MMA.  
\end{remark}

\begin{remark}\label{nonlocalremark}
A version of \ref{derivedMMCTforcAn} is actually true in the non-local case (see \cite{W}), but we explain here why the non-local case is a much more subtle problem.

For example, in the case $R=\mathbb{C}[u,v,x,y]/(uv-xy(y-2))$, define the modules $M_1:=(u,x)$, $L_1:=(u,x(y-2))$ and $L_2:=(u,xy)$.  Then we can show that both $M:=R\oplus M_1$ and $L:=R\oplus L_1\oplus L_2$ give NCCRs by checking complete locally. Also we can show  $\add M=\add L$ by checking complete locally \cite[2.26]{IW}). Thus $\End_R(M)$ and $\End_R(L)$ are Morita equivalent.  Since $\dim R=3$ and $R$ has an NCCR, by \cite{VdBNCCR} we obtain a derived equivalence between these NCCRs and some crepant resolution of $\Spec R$.   

Hence our picture in \ref{pictures} becomes
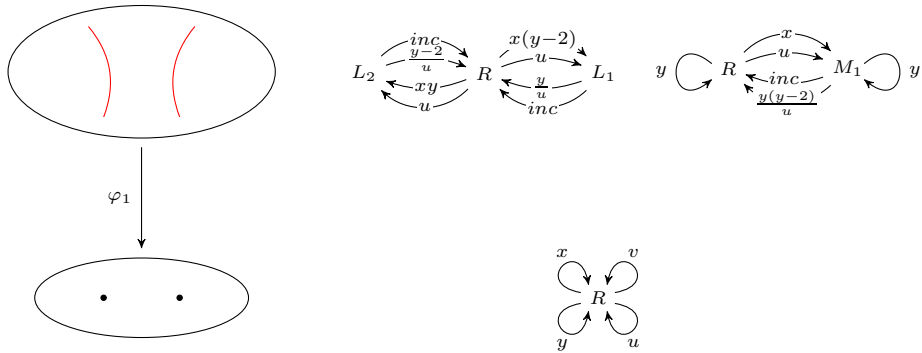
\begin{center}
\begin{tikzpicture}
\node (c2) at (0,3) {\begin{tikzpicture}
\draw (0.5,1.5) ellipse (50pt and 25pt);
\draw[red, bend left] (-0.2,2.1) to (0,0.9);
\draw[red, bend right] (1.2,2.1) to (1,0.9);
\end{tikzpicture}};

\node (d2) at (0,0)
{\begin{tikzpicture}
\draw (0,0) ellipse (40pt and 15pt);
\draw[fill] (0.5,0) circle (1pt);
\draw[fill] (-0.5,0) circle (1pt);
\end{tikzpicture}}; 

\node (b4) at (4.5,3)  
{\begin{tikzpicture}[xscale=0.9,yscale=0.9] 
\node (C1) at (0,0)  {$\scriptstyle R$};
\node (C1a) at (-0.1,0.05)  {};
\node (C1b) at (-0.1,-0.05)  {};
\node (C2) at (1.75,0)  {$\scriptstyle L_1$};
\node (C2a) at (1.85,0.05) {};
\node (C2b) at (1.85,-0.05) {};
\node (C3) at (-1.75,0)  {$\scriptstyle L_2$};
\node (C3a) at (-1.85,0.05) {};
\node (C3b) at (-1.85,-0.05) {};
\draw [->,bend left=45,looseness=1,pos=0.5] (C1) to node[inner sep=0.5pt,fill=white]  {$\scriptstyle x(y-2)$} (C2);
\draw [->,bend left=20,looseness=1,pos=0.5] (C1) to node[inner sep=0.5pt,fill=white]  {$\scriptstyle u$} (C2);
\draw [->,bend left=45,looseness=1,pos=0.5] (C2) to node[inner sep=0.5pt,fill=white]  {$\scriptstyle inc$} (C1);
\draw [->,bend left=20,looseness=1,pos=0.5] (C2) to node[inner sep=0.5pt,fill=white,below=-5pt]  {$\scriptstyle \frac{y}{u}$} (C1);
\draw [->,bend left=45,looseness=1,pos=0.5] (C1) to node[inner sep=0.5pt,fill=white]  {$\scriptstyle u$} (C3);
\draw [->,bend left=20,looseness=1,pos=0.5] (C1) to node[inner sep=0.5pt,fill=white]  {$\scriptstyle xy$} (C3);
\draw [->,bend left=45,looseness=1,pos=0.5] (C3) to node[inner sep=0.5pt,fill=white]  {$\scriptstyle inc$} (C1);
\draw [->,bend left=20,looseness=1,pos=0.5] (C3) to node[inner sep=0.5pt,fill=white,below=-5pt]  {$\scriptstyle \frac{y-2}{u}$} (C1);
\end{tikzpicture}};

\node (c4) at (8.5,3)  
{\begin{tikzpicture}[xscale=0.9,yscale=0.9] 
\node (C1) at (0,0)  {$\scriptstyle R$};
\node (C1a) at (-0.1,0.05)  {};
\node (C1b) at (-0.1,-0.05)  {};
\node (C2) at (1.75,0)  {$\scriptstyle M_1$};
\node (C2a) at (1.85,0.05) {};
\node (C2b) at (1.85,-0.05) {};
\draw [->,bend left=45,looseness=1,pos=0.5] (C1) to node[inner sep=0.5pt,fill=white]  {$\scriptstyle x$} (C2);
\draw [->,bend left=20,looseness=1,pos=0.5] (C1) to node[inner sep=0.5pt,fill=white]  {$\scriptstyle u$} (C2);
\draw [->,bend left=45,looseness=1,pos=0.5] (C2) to node[inner sep=0.5pt,fill=white]  {$\scriptstyle \frac{y(y-2)}{u}$} (C1);
\draw [->,bend left=20,looseness=1,pos=0.5] (C2) to node[inner sep=0.5pt,fill=white,below=-5pt]  {$\scriptstyle inc$} (C1);
\draw[<-]  (C1b) edge [in=140,out=-140,loop,looseness=10,pos=0.5] node[left] {$\scriptstyle y$} (C1b);
\draw[->]  (C2b) edge [in=-40,out=40,loop,looseness=10,pos=0.5] node[right] {$\scriptstyle y$} (C2b);
\end{tikzpicture}};

\node (d4) at (6,0)  
{\begin{tikzpicture} 
\node (1) at (0,0)  {$\scriptstyle R$};
\node (1a) at (-0.1,0.05)  {};
\node (1b) at (-0.1,-0.05)  {};
\node (2a) at (0.1,0.05)  {};
\node (2b) at (0.1,-0.05)  {};
\draw[<-]  (1b) edge [in=-170,out=-100,loop,looseness=12,pos=0.5] node[below] {$\scriptstyle y$} (1b);
\draw[<-]  (1a) edge [in=170,out=100,loop,looseness=12,pos=0.5] node[above] {$\scriptstyle x$} (1a);
\draw[->]  (2b) edge [in=-80,out=-10,loop,looseness=12,pos=0.5] node[below] {$\scriptstyle u$} (2b);
\draw[->]  (2a) edge [in=80,out=10,loop,looseness=12,pos=0.5] node[above] {$\scriptstyle v$} (2a);
\end{tikzpicture}};

\draw[->] (0,2) -- node[left] {$\scriptstyle \varphi_1$} (d2);

\end{tikzpicture}
\end{center}
with everything in the top row being derived equivalent.   This all relies on the fact that we a priori know that we have an NCCR (because we can check this complete locally) and so we can use \cite{VdBNCCR}.  In a similar case $uv=(x^2+y^3)y(y-2)$ (where no NCCR exists), we cannot play this trick, since we do not know if MMAs can be detected locally, and also since completion destroys $\mathds{Q}$-factoriality we cannot reduce to the complete local setting. 

The other way to try to prove a non-local version of \ref{derivedMMCTforcAn} is to try and find an explicit tilting complex starting with the geometry, as in \ref{Dbforflags}.  Although a tilting bundle exists abstractly for one-dimensional fibres \cite{VdB1d}, to describe it explicitly is delicate, since for example we have to deal with issues such as curves being forced to flop together.   We refer the reader to \cite{W} for more details.
\end{remark}

\end{document}